\documentclass[12pt,reqno]{amsart}

\usepackage[margin=1in]{geometry}
\usepackage{mathtools}
\usepackage{amssymb}
\usepackage{amsthm}
\usepackage{dsfont}
\usepackage{graphicx} 
\usepackage{microtype}

\usepackage{tikz}
\usetikzlibrary{shapes.geometric, positioning, fit}

\usepackage{charter}
\usepackage[cal=euler]{mathalpha}

\usepackage{xcolor}
\usepackage[backref=page]{hyperref}
\hypersetup{
    colorlinks=true,
    linkcolor=cyan!80!black,
    citecolor=cyan!80!black,
    urlcolor=magenta,
}

\usepackage{booktabs}
\usepackage{makecell}
\usepackage[symbol]{footmisc}

\usepackage{apptools}
\AtAppendix{\counterwithin{lemma}{section}}

\allowdisplaybreaks

\newcommand{\bbE}{\mathbb{E}}

\newcommand{\bbN}{\mathbb{N}}
\newcommand{\bbP}{\mathbb{P}}
\newcommand{\bbR}{\mathbb{R}}

\newcommand{\bbI}{\mathbb{I}}

\newcommand{\cT}{\mathcal{T}}

\newcommand{\cW}{\mathcal{W}}

\newcommand{\cG}{\mathcal{G}}
\newcommand{\cE}{\mathcal{E}}

\newcommand{\cV}{\mathcal{V}}

\newcommand{\Binom}{\mathsf{Binomial}}

\newcommand{\Unif}{\mathsf{Uniform}}
\newcommand{\Ber}{\mathsf{Bernoulli}}

\newcommand{\Poisson}{\mathsf{Poisson}}

\newcommand{\ER}{\mathrm{ER}}
\newcommand{\dist}{\mathrm{dist}}
\newcommand{\iid}{\mathrm{i.i.d.}}
\newcommand{\cc}{\mathrm{CC}}
\newcommand{\PGWT}{\mathrm{PGWT}}
\newcommand{\PPT}{\mathrm{PPT}}
\newcommand{\PPP}{\mathrm{PPP}}

\DeclareMathOperator{\Var}{Var}

\DeclareMathOperator{\Cov}{Cov}

\newcommand{\convp}{\xrightarrow{\bbP}}
\newcommand{\convd}{\xrightarrow{d}}

\newcommand{\eqas}{\overset{\text{a.s.}}{=}}

\newcommand{\lf}{\left(}
\newcommand{\rf}{\right)}

\newcommand{\lt}{\left[}
\newcommand{\rt}{\right]}

\let\hat\widehat
\let\tilde\widetilde

\theoremstyle{plain}
\newtheorem{theorem}{Theorem}[section]
\newtheorem{corollary}{Corollary}[section]

\newtheorem{lemma}{Lemma}[section]

\theoremstyle{definition}
\newtheorem{assumption}{Assumption}[section]
\newtheorem{definition}{Definition}[section]

\theoremstyle{remark}
\newtheorem{remark}{Remark}[section]

\title[Local convergence of Triadic Closures ]{Local convergence of Triadic Closures of Sparse Random Graphs}
\author[H. Talukdar]{Himasish Talukdar}
\address{
    Statistics and Mathematics Unit \\
    Indian Statistical Institute \\
    203 B.T. Road, Kolkata 700108 \\
    West Bengal, India.
}
\email{talukdar.himasish@gmail.com}

\begin{document}

\begin{abstract}
Triadic closure is a natural mechanism for modeling the clustering phenomenon observed in many real-world networks. Given a possibly random undirected graph $H$, we consider pairs of non-adjacent vertices that have a common neighbor. Conditional on $H$, each such pair is joined by an edge independently and simultaneously, with a probability depending on the degrees of its endpoints and on the number of their common neighbors. The resulting graph is called a triadic closure of $H$. We study the effect of this operation on local convergence. More precisely, we prove that if a sequence of potentially random finite undirected graphs of $\left( H_n \right)_{n \in \bbN}$ converges locally in probability to a rooted graph $H^*$, then the corresponding sequence of triadic closures converges locally in probability to the same triadic closure of $H^*$. As a consequence, we obtain a general limit theorem for the clustering coefficients of triadic closures. In particular, we show that the triadic closures can have a strictly positive limiting clustering coefficients even when the limiting clustering coefficients of the parent graph sequence is zero. Finally, we apply our results to several classical sparse random graph models including Erd\H{o}s-R\'{e}nyi, configuration, preferential attachment models and compute the limiting clustering coefficients of their triadic closures.

\end{abstract}

\keywords{Local convergence of random graphs; Clustering coefficients; Sparse random graphs}
\subjclass[2020]{05C80, 05C82, 60B99}

\maketitle
\thispagestyle{empty}
\tableofcontents

\section{Introduction}\label{sec:intro}

High clustering is a characteristic feature of many real networks. More precisely, in real networks it is often the case that the neighbors of a node are connected to one another. For example, \cite{watts1998collective} reported local clustering coefficient of $0.79$ for a film-actor collaboration network and $0.28$ for the neuronal network of Caenorhabditis elegans. Classical sparse models, such as Erd\H{o}s--R\'{e}nyi and random regular graphs, have asymptotically vanishing clustering and therefore fail to reproduce this feature. A natural way to introduce clustering is \emph{triadic closure}, whereby vertices sharing a neighbor form a new connection, completing a triangle. This mechanism reflects the formation of relationships through mutual acquaintances and is supported by observations of scientific collaborations~\cite{newman2001clustering}. We formalize it through a probabilistic closure operation on a parent graph. Clustering coefficients alone, however, do not describe how triangles and other cycles overlap within neighborhoods. We therefore study local convergence of the resulting graphs, which captures the distribution of neighborhoods around a uniformly chosen vertex, and derive limits for local and global clustering coefficients as a consequence.

\subsection{Our model}
Let $H$ be a finite, possibly random, simple graph, called the \emph{parent graph}, and let $\phi:\bbN^3\to[0,1]$ be symmetric in its first two coordinates. Conditional on $H$, each unordered pair of distinct, nonadjacent vertices $u,v$ with a common neighbor receives an edge with probability
\[
    \phi\left(d_H(u),d_H(v),\tilde{d}_H(\{u,v\})\right),
\]
independently of the other pairs. Here $d_H(u)$ is the degree of $u$, and $\tilde{d}_H(\{u,v\})$ is the number of common neighbors of $u$ and $v$. All edges of $H$ are retained. The additions are simultaneous, with probabilities determined entirely by $H$, and the resulting graph is denoted by $\cT_\phi(H)$. Thus, closure may depend on both endpoint degrees and the number of common neighbors. We write $\cT_p$ when $\phi\equiv p$, i.e., the closure probability is constant.

\subsection{Related works}
Triadic closure has been incorporated into several models of network growth. \cite{holme2002growing} added triad formation to preferential attachment to obtain tunable clustering. \cite{vazquez2003growing} studied growth mechanisms based on local information, while \cite{jackson2007meeting} combined random meetings (i.e. edges) with meetings through existing neighbors. In these models, closure is part of the growth dynamics. Static constructions provide another approach. For example, \cite{newman2009random} used individual edges and triangles as building blocks, while random intersection models connect vertices that share an attribute or belong to the same group~\cite{deijfen2009random,bloznelis2013degree}.

More directly related are models that apply closure to a sampled parent graph. \cite{zhang2023generalized} combined a generalized configuration model with random triadic closure, allowing community structure and degree correlations. \cite{cirigliano2024strongly} introduced a static closure model and derived degree distributions, global clustering and small-subgraph counts for uncorrelated, locally tree-like backbones (i.e. parent graphs). For the same class of backbones, \cite{cirigliano2026strongly} subsequently obtained degree correlations and the clustering spectrum, whose average over the final degree distribution gives the local clustering coefficient. Their closure operation corresponds to $\phi(a,b,c)=1-(1-p)^c$. On a tree, it agrees with $\cT_p$, since every admissible pair has exactly one common neighbor.

Our approach uses local convergence, which describes the asymptotic distribution of neighborhoods around a uniformly chosen vertex. Foundational works in this direction include \cite{benjamini2001recurrence,aldous2004objective,aldous2007processes}. For a detailed study of local convergence of random graphs see \cite{banerjee2024local,vanderhofstad2024random}. The connection between local limits and clustering is already part of this theory. This connection has been exploited in many earlier works including \cite{vanderhofstad2017local}, \cite{vanderhofstad2018triadic} and \cite{vanderhofstad2023scaling} where the authors considered the hidden-variable models, erased configuration models, and spatial inhomogeneous random graphs respectively. \cite{kurauskas2022local} related local convergence, degree moments and subgraph counts, giving criteria under which graph statistics are determined by the local limit. In this paper, we also use this connection to derive limiting formulae for local and global clustering coefficients from the local convergence of the triadic closure graphs.

\subsection{Our contributions}
Theorem~\ref{thm:lwl} shows that if $|\cV(H_n)|\to\infty$ and $H_n$ converges locally in probability to a locally finite rooted graph $H^*\sim\mu$, then
\[
    G_n:=\cT_\phi(H_n)
    \quad\text{converges locally in probability to}\quad
    G^*:=\cT_\phi(H^*).
\]
The law $\mu$ may be random. The theorem allows cycles and dependence between neighboring degrees in the parent limit, and requires no degree-moment assumptions. It also extends to any fixed finite number of closure steps, with possibly different weight functions (Corollary~\ref{cor:lwl_iterative}). The proof uses the fact that the conditional law of a radius-$r$ neighborhood after closure is determined by the radius-$(2r+1)$ parent neighborhood, together with control of the fluctuations of empirical neighborhood frequencies.

Theorem~\ref{thm:cc_local} then gives convergence in probability of the local clustering coefficient and an explicit expression for its limit in terms of $H^*$ and $\phi$. Corollary~\ref{cor:degree_distribution} identifies the limiting degree distribution. For global clustering coefficient, under some extra degree moment assumption, we prove convergence in probability and express the limit through expected triangle and wedge counts in $G^*$ (see Theorem~\ref{thm:cc_global}). The formulae simplify for tree limits and constant closure probabilities (see Corollary~\ref{cor:cc_local} and \ref{cor:cc_global}).

For $\phi\equiv p$, we apply the results to sparse Erd\H{o}s--R\'{e}nyi graphs, random regular graphs, configuration models and preferential attachment graphs. For the first three families, the global clustering formulae agree with those of~\cite{cirigliano2024strongly}, and our results establish convergence in probability under the stated hypotheses. The preferential attachment example uses the P\'{o}lya-point tree and accounts for dependence arising from vertex ages and attachment histories. These examples show how closure can produce positive limiting clustering from locally tree-like parent graphs when $p>0$. Our central contribution is therefore a local convergence theorem for general weight function, together with explicit clustering limits that also apply to cyclic or degree correlated parent limits.

\subsection{Notational conventions}

Let $\bbN$ denote the set of natural numbers. For $n \in \bbN$, let $[n]$ denote the set $\{1, 2, \ldots, n\}$. For a set $A$ and positive integer $k$, $\binom{A}{k}$ denotes the collection of $k$-subsets of $A$. So, in context of an undirected graph with vertex set $V$, the collection of all possible edges is given by $\binom{V}{2}$. For a graph $H$, let $\cV(H)$ and $\cE(H)$ denote the set of vertices and edges of $H$ respectively. For vertices $u, v \in \cV(H)$, let $d_H(u)$ denote the degree of $u$ in $H$ and $\tilde{d}_{H}(u,v)$ denote the number of common neighbors of $u$ and $v$ in $H$. Moreover, for two vertices $u, v \in \cV(H)$, let $\dist_H (u, v)$ denote the graph distance between $u$ and $v$ in $H$. For $v \in \cV(H)$ and $r >0$, let $B_H (v, r)$ denote the closed ball of radius $r$ around the point $v$ in the graph $H$, i.e.,
\begin{align*}
    \cV(B_H (v, r)) &:= \{u \in \cV(H) : \dist_H (u, v) \le r\}, \\
    \cE(B_H (v, r)) &:= \{ \{u, w\} \in \cE(H): u, w \in \cV (B_H(v,r)) \}.
\end{align*}
Define the boundary of $B_H (v, r)$, denoted by $\partial B_H (v, r)$, as 
\begin{align*}
    \cV (\partial B_H (v, r)) &:= \{u \in \cV(H) : \dist_H (u, v) = r\}, \\
    \cE(\partial B_H (v, r)) &:= \{ \{u, w\} \in \cE(H): u, w \in \cV (\partial B_H(v,r)) \}.
\end{align*}
A graph of the form $(\{u,v,w\}, \{\{u,v\}, \{v,w\}\})$ is called a \emph{wedge}. The \emph{wedge set of $H$}, denoted by $\cW(H)$ is the collection of edges $\{u, v\} \in \binom{\cV(H)}{2}$ such that $\{u, v\} \notin \cE(H)$, but there exists $w \in \cV(H)$ with $\{ u, w\}, \{w, v\} \in \cE(H)$. The reader should note that the elements of the wedge set, as per our definition, are not wedges themselves, rather they are edges.

For two graphs $H_1, H_2$, by $H_1 \cong H_2$, we mean that $H_1$ and $H_2$ are isomorphic as graphs. More precisely, $H_1 \cong H_2$ if there exists a bijection $\sigma : \cV(H_1) \to \cV(H_2)$ such that $\{u, v\} \in \cE(H_1)$ if and only if $\{\sigma(u), \sigma(v)\} \in \cE(H_2)$. A rooted graph $H^*$ consists of a vertex set $\cV(H^*)$, an edge set $\cE(H^*)$ and a prespecified vertex $o_{H^*}$, called the root. In this text, the notations used for rooted graphs use a $*$ in the superscript to distinguish it from graphs. Two rooted graphs $H_1^*$, $H_2^*$ are called isomorphic, denoted by $H_1^* \cong H_2^*$, if there exists a bijection $\sigma : \cV(H_1^*) \to \cV(H_2^*)$ such that $\sigma(o_{H_1^*}) = o_{H_2^*}$ and $\{u,v\} \in \cE(H_1^*)$ if and only if $\{\sigma(u), \sigma(v)\} \in \cE(H_2^*)$.

Throughout the text, $\frac{0}{0}$ is interpreted as $0$, an empty sum is interpreted as $0$ and an empty product is interpreted as $1$.

\subsection{Organisation of the paper}
Section~\ref{sec:model} introduces the triadic closure model in detail. Section~\ref{sec:prelim} defines local convergence of random graphs and the two clustering coefficients. The main results are stated in Section~\ref{sec:main_results} and proved in Section~\ref{sec:proofs}. Section~\ref{sec:empirical} contains an empirical study. Section~\ref{sec:examples} applies the results to Erd\H{o}s-R\'{e}nyi graphs, random regular graphs, configuration models and preferential attachment graphs. Finally, Appendix~\ref{sec:aux} contains the proofs of a few supporting lemmas.

\section{The model}\label{sec:model}

In this section, we describe the model in detail. By a \textit{weight function} we mean a function $\phi: \bbN^3 \to [0,1]$ which is symmetric in the first two coordinates. Given a graph $H$ and a weight function $\phi$, define $\phi_{H}: \cW(H) \to [0,1]$, by 
\[
    \phi_H(e) := \phi(d_{H}(u), d_{H}(v), \tilde{d}_H(\{u,v\})).
\]
Recall that $\tilde{d}_{H}(\{u,v\})$ denotes the number of common neighbors of $u$ and $v$ in $H$.  Note that $\phi_H$ is well defined since for any $\{u,v\} \in \cW(H)$, $d_{H}(u), d_{H}(v), \tilde{d}_H(u,v) \ge 1$. We start with the definition of triadic closures of graphs.
\begin{definition}[$\phi$-triadic closure]
    Let $H$ be a (possibly random) graph and $\phi$ be a weight function. The \emph{$\phi$-triadic closure of $H$}, denoted by $\cT_\phi(H)$, is a random graph with 
    \begin{align*}
        &\cV(\cT_\phi(H)) := \cV(H),
    \end{align*}
    and for $e := \{ u, v\} \in \binom{\cV(H)}{2}$,
    \begin{align*}
        e &\in \cE(\cT_\phi(H)) \qquad\,\, \text{if} \,\, e \in H,\\
        e &\in \cE(\cT_\phi(H)) \qquad \text{with probability} \,\, \phi_H(e)\,\, \text{if} \,\, e \in \cW(H),\\
        e &\notin \cE(\cT_\phi(H)) \qquad \text{with probability} \,\, 1 - \phi_H(e)\,\, \text{if} \,\, e \in \cW(H),\\
        e &\notin \cE(\cT_\phi(H)) \qquad \text{if} \,\, \dist_H(u, v) \ge 3.
    \end{align*}
\end{definition}
In other words, conditional on $H$, define a collection of independent random variables $Z_H:= (Z_H(e))_{e \in \cW(H)}$ with $Z_H (e) \sim \Ber(\phi_H(e))$ and define 
\[
    \bbI(e \in \cT_\phi(H)) := Z_H (e) \qquad \text{for}\,\, e \in \cW(H).
\]
Note that $\cT_\phi(H)$ is a supergraph of $H$ on the same vertex set. Essentially, each wedge is closed with a probability depending on the degrees of its endpoints (i.e. $u$ and $w$) and everything else is kept fixed. When $\phi$ is clear from the context, we just write $\cT(H)$ for $\cT_\phi(H)$. If $\phi \equiv p$ for some $p \in [0,1]$, $\cT_\phi$ is denoted by $\cT_p$. We call $\cT_1(H)$ the \emph{complete triadic closure} of $H$. On the other hand, notice that $\cT_0(H) = H$. For a rooted graph $H^* := (\cV(H^*), \cE(H^*), o_{H^*})$, the $\phi$-triadic closure of $H^*$, $\cT_\phi(H^*)$, is a random rooted graph given by
\[
    \cT_\phi(H^*) := (\cT_\phi (\cV(H^*), \cE(H^*)), o_{H^*}).
\]
So, the $\phi$-triadic closure of a rooted graph is the $\phi$-triadic closure of the corresponding unrooted graph with the same root. 

In this article we explore the geometric properties of triadic closures of random graphs. The necessary concepts in this respect are introduced in the next section.

\section{Local convergence and clustering coefficients}\label{sec:prelim}

In this section we define the concept of local convergence and clustering coefficients which appear in the main results of this article. The concepts defined below naturally generalises to multigraph setting. However, for simplicity we only explain these concepts for graphs. The reader is referred to \cite{vanderhofstad2024random} for a comprehensive discussion on this. 

\subsection{Local convergence}

The concept of local convergence (also known as Benjamini- Schramm convergence) of graphs was first introduced by Benjamini and Schramm in their seminal work \cite{benjamini2001recurrence}. Since then there has been numerous works on it \cite{aldous2004objective, aldous2007processes, bordenave2010resolvent, bordenave2015large, kurauskas2022local, banerjee2024local, lo2024local}. Here we give a brief description of local convergence of graphs.

We have already defined rooted graphs. Our first goal is to introduce a metric on the space of rooted graphs. In order to do so, we require an equivalence on the space of rooted graphs.For a rooted graph $H^*$, let $\mathcal{C}_{H^*}(o_{H^*})$ denote the connected component of $o_{H^*}$ in $H^*$ seen as a rooted graph with the same root $o_{H^*}$. Now two rooted graphs $H_1^*$ and $H_2^*$ are equivalent if $\mathcal{C}_{H_1^*}(o_{H_1^*})$ and $\mathcal{C}_{H_2^*}(o_{H_2^*})$ are isomorphic as rooted graphs. Let $\cG^*$ denote the set of equivalence classes of rooted locally finite graphs under the mentioned equivalence relation. Note that each equivalence class in $\cG^*$ contains exactly one rooted connected graph up to isomorphism. To make things concrete we take this rooted graph to be the representative of each class. Now, we define a metric $\rho$ on $\cG^*$ as follows:
For $H_1^*, H_2^* \in \cG^*$, 
\[
    \rho(H_1^*, H_2^*) := \frac{1}{1 + R(H_1^*, H_2^*)},
\]
where
\[
    R(H_1^*, H_2^*) := \sup \{ r >0 : B_{H_1^*}(o_{H_1^*}, r) \cong B_{H_2^*}(o_{H_2^*}, r)\}.
\]
One can show that $(\cG^*, \rho)$ is a Polish space (see \cite{vanderhofstad2024random}). Next we define local convergence of random graphs. Suppose $(H_n)_{n \in \bbN}$ is a sequence of finite, possibly random graphs. 
\begin{definition}[Local weak convergence]\label{def:lwc}
    \begin{enumerate}
        \item[(a)] We say that $H_n$ converges \emph{locally weakly} to $H^* \sim \mu$, where $\mu$ is a deterministic measure on $\cG^*$, if for every bounded and continuous $f: \cG^* \to \bbR$, we have
        \[
            \bbE \left[ \frac{1}{|\cV(H_n)|} \sum_{v \in \cV(H_n)}f((H_n, v)) \right] \to \bbE_\mu [f(H^*)].
        \]
        \item[(b)] We say that $H_n$ converges \emph{locally in probability} to $H^* \sim \mu$, where $\mu$ is a possibly random measure on $\cG^*$, if for every bounded and continuous $f: \cG^* \to \bbR$, we have
        \[
            \frac{1}{|\cV(H_n)|} \sum_{v \in \cV(H_n)}f((H_n, v)) \convp \bbE_\mu [f(H^*)].
        \]
    \end{enumerate}    
\end{definition}
In either of the definitions the limiting measure $H^*$ is called the \emph{local limit} of the sequence $H_n^*$. Moreover, it is not hard to see that definition~\ref{def:lwc}(b) is a stronger notion of convergence than definition~\ref{def:lwc}(a) if $\mu$ is a deterministic measure. More generally, if $H_n$ converges locally in probability to $H^* \sim \mu$, where $\mu$ is a possibly random measure on $\cG^*$, then $H_n$ converges locally weakly to $\overline{H}^* \sim \overline{\mu}$, where the deterministic measure $\overline{\mu}$ is defined by its action on the test functions as follows:
\begin{equation}\label{eq:prob_weak}
    \bbE_{\overline{\mu}} [f(\overline{H}^*)] := \bbE[\bbE_{\mu} [f(H^*)]],
\end{equation}
for all $f \in C_b(\cG^*, \rho)$. 
\begin{remark}
    If $\mu$ is a random measure, then $\bbE_\mu [f(H^*)] $ is a random variable and for each fixed $\omega$ in the underlying probability space $\bbE_{\mu(\omega)} [f(H^*)]$ is a real number. More precisely
    \begin{align*}
        \bbE_{\mu(\omega)} [f(H^*)] = \int_{\cG^*} f(H^*) \mu(\omega) (d H^*)
    \end{align*}
    So, the outer expectation in the RHS of \eqref{eq:prob_weak} averages out the randomness in $\mu$ over the variable $\omega$.
\end{remark}

We now give an equivalent characterization of local convergence. To that end, we define subgraph density.

\begin{definition}[Subgraph density]\label{def:subgraph_density}
    Given a finite graph $H$, $K^* \in \cG^*$ and $r \ge 0$, the \emph{subgraph density} at radius $r$ of $K^*$ in $H$ is defined as
    \begin{equation}
        p_H(K^*, r) := \frac{1}{|\cV(H)|} \sum_{v \in \cV(H)} \bbI ((B_H(v, r), v) \cong K^*).
    \end{equation}
\end{definition}

Now we mention equivalent criteria of local convergence to be used in the proof of Theorem~\ref{thm:lwl}.

\begin{theorem}[Theorem~2.15 of \cite{vanderhofstad2024random}]\label{thm:lwc_equivalent_criterion}
    \begin{enumerate}
        \item[(a)] $H_n$ converges locally weakly to $H^* \sim \mu$, where $\mu$ is deterministic measure, if and only if for all $K^* \in \cG^*$ and $r \ge 0$,
        \begin{equation}\label{eq:lwc_equivalent}
            \bbE \lt p_{H_n}(K^*, r) \rt \to \mu(B_{H^*}(o_{H^*}, r) \cong K^*).
        \end{equation}
        \item[(b)] $H_n$ converges locally in probability to $H^* \sim \mu$, where $\mu$ is a possibly random measure, if and only if for all $K^* \in \cG^*$ and $r \ge 0$,
        \begin{equation}\label{eq:lwc_in_prob_equivalent}
            p_{H_n}(K^*, r) \convp \mu(B_{H^*}(o_{H^*}, r) \cong K^*).
        \end{equation}
    \end{enumerate}
\end{theorem}

By virtue of Theorem~\ref{thm:lwc_equivalent_criterion}, roughly, local convergence means that the neighborhoods of a randomly chosen vertex in $H_n$ looks like the corresponding neighborhoods of $o_{H^*}$ in $H^*$.
\subsection{Clustering coefficients}
In this subsection, we introduce the concept of clustering coefficients and provide its connection with local convergence. We start with the definitions of global and local clustering coefficients.

Roughly speaking, clustering coefficients measure in a graph, on average, what proportion of wedges are closed to form triangles. Or, in the context of social networks, it can be thought of as what is the probability that two random friends of a random individual are friends themselves. For example, since a tree does not have any triangle, the clustering coefficients should be $0$. On another extreme, the clustering coefficients of a complete graph (with at least three vertices) should be $1$. We now formally define clustering coefficients of a graph.

Suppose $H$ is a finite graph. In what follows, $\frac{0}{0}$ is interpreted as $0$. Denote by $W_{H}$, twice the number of wedges (incomplete or complete) in $H$, i.e.,
\begin{equation}\label{eq:wedge_count_defn}
    W_H:= \sum_{\substack{u,v,w \in \cV(H)\\ u \ne w}} \bbI(\{u,v\}, \{v, w\} \in \cE(H)),
\end{equation}
which further simplifies to 
\[
    W_H = \sum_{v \in \cV(H)} d_H(v) (d_H(v) - 1) = \sum_{v \in \cV(H)} d_H(v)^2 -2 |\cE(H)|.
\]
In the same spirit, denote by $T_{H}$, six times the number of triangles in $H$, i.e.,
\[
    T_H := \sum_{u,v,w \in \cV(H)} \bbI(\{u,v\}, \{v, w\}, \{w, u\} \in \cE(H)).
\]
\begin{definition}\label{def:cc_global}
    The \emph{global clustering coefficient of $H$}, denoted by $\cc_H$, is defined as
    \[
        \cc_H := \frac{T_H}{W_H}.
    \]
\end{definition}
Next, we define the concept of local clustering coefficient. Let $v \in \cV(H)$. Denote by $T_H(v)$, twice the number of triangles $v$ is a part of, or formally,
\[
    T_H(v) := \sum_{u, w \in \cV(H)} \bbI(\{u, v\}, \{v, w\}, \{w, u\} \in \cE(H)).
\]
Then define the local clustering coefficient for $v$ to be
\[
    \cc_H(v) := \frac{T_H(v)}{d_H(v) (d_H(v) - 1)}.
\]
\begin{definition}\label{def:cc_local}
    The \emph{local clustering coefficient of $H$}, denoted by $\overline{\cc}_H$, is defined as  
    \[
        \overline{\cc}_H := \frac{1}{|\cV(H)|} \sum_{v \in \cV(H)} \cc_H(v).
    \]
\end{definition}
On the other hand if $H^* \sim \mu$ is a possibly random, possibly infinite rooted graph, then these concepts are defined as follows.
\begin{definition}\label{def:cc_global_local_rooted}
The global and local clustering coefficients of $H^*$ are defined as
\begin{align*}
    \cc_{H^*} &:= \frac{\bbE_\mu [T_{H^*}(o^*)]}{\bbE_\mu [d_{H^*}(o_{H^*}) (d_{H^*}(o_{H^*}) - 1)]}, \\
    \overline{\cc}_{H^*} &:= \bbE_\mu \lt\frac{ T_{H^*}(o^*)}{ d_{H^*}(o_{H^*}) (d_{H^*}(o_{H^*}) - 1)} \rt,
\end{align*}
respectively, whenever the RHS is well defined.
\end{definition}
The following theorems from \cite{vanderhofstad2024random} establishes the relationship between local convergence and the convergence of clustering coefficients. 
\begin{theorem}[Theorem~2.22 of \cite{vanderhofstad2024random}]\label{thm:lwc_cc_global_remco}
    Let $(H_n)_{n \in \bbN}$ be sequence of finite graphs such that $|\cV(H_n)| \to \infty$. Assume that $H_n$ converges locally in probability to
    $H^* \sim \mu$, where $\mu$ is a possibly random measure on $\cG^*$. Further assume that $\left(d_{H_n}(o_{H_n})^2 \right)_{n \in \bbN}$ is uniformly integrable where $o_{H_n} \sim \Unif(\cV(H_n))$ and that $\mu(d_{H^*}(o_{H^*}) >1) >0$. Then,
    \begin{equation}\label{eq:lwc_cc_global_remco}
        \cc_{H_n} \convp \frac{\bbE_\mu [ T_{H^*}(o_{H^*}) ]}{\bbE_\mu [d_{H^*}(o_{H^*})(d_{H^*}(o_{H^*}) - 1)]}. 
    \end{equation}
\end{theorem}
\begin{theorem}[Theorem~2.23 of \cite{vanderhofstad2024random}]\label{thm:lwc_cc_local_remco}
    Let $(H_n)_{n \in \bbN}$ be sequence of graphs such that $|\cV(H_n)| \to \infty$. Assume that $H_n$ converges locally in probability to $H^* \sim \mu$, where $\mu$ is a possibly random measure on $\cG^*$. Then,
    \begin{equation}\label{eq:lwc_cc_local_remco}
        \overline{\cc}_{H_n} \convp \bbE_\mu \left[ \frac{T_{H^*}(o_{H^*}) }{d_{H^*}(o_{H^*})(d_{H^*}(o_{H^*}) - 1)} \right]. 
    \end{equation}
\end{theorem}
Now we are ready to state our main results.

\section{Main results}\label{sec:main_results}
\subsection{Local convergence}
Our first theorem essentially says that $\phi$-triadic closures act continuously under local convergence. More precisely,

\begin{theorem}[Local limit]\label{thm:lwl}
    Suppose $(H_n)_{n \in \bbN}$ is a sequence of finite graphs such that $|\cV(H_n)| \to \infty$. Assume that $H_n$ converges locally in probability to $H^* \sim \mu$, where $\mu$ is a possibly random measure on $\cG^*$. Define $G_n := \cT_\phi (H_n)$ for some weight function $\phi$. Then, $G_n$ converges locally in probability to $\cT_\phi (H^*)$.
\end{theorem}

As an immediate corollary of Theorem~\ref{thm:lwl}, we get the following.

\begin{corollary}\label{cor:lwl_iterative}
    If $\phi_1, \ldots, \phi_k$ are $k$ weight functions, then under the assumption of Theorem~\ref{thm:lwl}, we have that $\cT_{\phi_k} \circ \cdots \circ \cT_{\phi_1} (H_n)$ converges locally in probability to $\cT_{\phi_k} \circ \cdots \circ \cT_{\phi_1} (H^*)$.
\end{corollary}

\begin{remark}\label{rmk:lwl}
    As the proof of Theorem~\ref{thm:lwl} will show, similar statements hold for other graphs operations instead of triadic closure as long as this operation affects only a fixed finite neighborhood of the parent graph. For example, instead of taking the closure probabilities to only depend on the degrees and the number of common neighbors of the end points, one can make it a function of number of certain fixed subgraph (e.g. triangles) formed at the end points in the parent graph. 
\end{remark}

\subsection{Clustering coefficients}\label{subsec:main_results_cc}

Before we dive into the main results of this subsection, we introduce a few definitions and notations which will make the statements rather concise.

\begin{definition}[Induced wedge set]\label{def:induced_wedge_set}
    Let $G^*$ be a rooted graph and $U, W \subseteq \cV(G^*)$. Then, $\cW(G^*, \{U, W\})$, the \emph{wedge set induced by $\{U, W\}$ in $G^*$} is the collection of edges $e \in \cW(G^*)$ such that one endpoint of $e$ is in $U$ and the other is in $W$. In case $U=W$, we just write $\cW(G^*, U)$.
\end{definition}

Let $G^*:= \cT_\phi(H^*)$. For $v \in \cV(\partial B_{H^*}(o_{H^*}, 2))$, define 
\[
    X_v:= Z_{H^*}(\{o_{H^*}, v\})= \bbI( \{o_{H^*}, v\} \in \cE(G^*)).
\]
Then, given $H^*$, 
\[
    X_v \sim \Ber(\phi_{H^*}(\{o_{H^*}, v\})),
\]
and are independent for $v \in \cV(\partial B_{H^*}(o_{H^*}, 2))$. Also, define
\[
    X:= \sum_{v \in \cV(\partial B_{H^*} (o_{H^*}, 2))} X_v.
\]
Note that there are $\frac{1}{2} T_{H^*}(o_{H^*})$ triangles containing the root $o_{H^*}$ in $H^*$. Note that the corresponding set of edges is given by 
\[
    \{ \{u,v\} \in \cE(H^*): u, v \in \partial B_{H^*}(o_{H^*}, 1) \}.
\]
Define
\[
    F_1 := \binom{\cV(\partial B_{H^*}(o_{H^*}, 1))}{2} \setminus \{ \{u,v\} \in \cE(H^*): u, v \in \partial B_{H^*}(o_{H^*}, 1) \}.
\]
Alternatively, in the notation introduced in Definition~\ref{def:induced_wedge_set},\[
    F_1 = \cW (B_{H^*}(o_{H^*}, 1), \cV(\partial B_{H^*}(o_{H^*}, 1))),
\]
the wedge set in $B_{H^*}(o_{H^*}, 1)$ induced by $\cV(\partial B_{H^*}(o_{H^*}, 1))$. We also define
\begin{align*}
    F_2 &:= \cW (H^*, \{ \cV(\partial B_{H^*}(o_{H^*}, 1)), \cV(\partial B_{H^*}(o_{H^*}, 2))\}) , \\
    F_3 &:= \cW (H^*, \cV(\partial B_{H^*}(o_{H^*}, 2))).
\end{align*}
Now, conditional on $H^*$, $Z_{H^*}(e) \sim \Ber(\phi_{H^*}(e))$ for $e \in F_1 \cup F_2 \cup F_3$ and are independent for $e \in F_1 \cup F_2 \cup F_3$. 

Lastly, conditional on $H^*$, the collections $(X_v)_{v \in \partial B_{H^*}(o_{H^*}, 2)}$ and $(Z_{H^*}(e))_{e \in F_1 \cup F_2 \cup F_3}$ are independent. 
We are now ready to state the main theorem of this subsection which describes the limit of the local clustering coefficients. 

\begin{theorem}[Local clustering coefficients]\label{thm:cc_local}
    Let $(H_n)_{n \in \bbN}$ be sequence of finite graphs such that $|\cV(H_n)| \to \infty$. Assume that $H_n$ converges locally in probability to $H^* \sim \mu$, where $\mu$ is a possibly random measure on $\cG^*$. Let $\phi$ be a weight function and $G_n:= \cT_\phi(H_n)$ for $n \in \bbN$ and $G^*:= \cT_\phi(H^*)$. Then,
    \begin{equation}\label{eq:cc_local}
        \overline{\cc}_{G_n} \convp \bbE_\mu \left[ \frac{T}{D} \right],
    \end{equation}
    where
    \begin{align}
        X & := \sum_{v \in \cV(\partial B_{H^*} (o_{H^*}, 2))} X_v, \label{eq:X_defn}\\
        D &:= d_{H^*}(o_{H^*}) \cdot (d_{H^*}(o_{H^*}) - 1)  + (2 d_{H^*}(o_{H^*}) -1)\cdot X + X^2, \label{eq:D_defn}\\
         T &:= T_{H^*}(o_{H^*}) + 2 \sum_{e \in F_1} \phi_{H^*}(e) + 2 \sum_{v \in \cV(\partial B_{H^*}(o_{H^*}, 2))} \left( d_{B_{H^*}(o_{H^*}, 2)} (v) - d_{\partial B_{H^*}(o_{H^*}, 2)} (v) \right) \cdot X_v \label{eq:T_defn}\\
         &+ 2 \sum_{\substack{ \{u,v\} \in F_2 \\ v \in \cV(\partial B_{H^*}(o_{H^*}, 2))}} \phi_{H^*}(\{u,v\}) \cdot X_v + \sum_{\substack{u, v \in \cV(\partial B_{H^*}(o_{H^*}, 2))\\ u \ne v}} \bbI(\{u, v\} \in \cE(H^*)) \cdot X_u X_v \nonumber\\
         &+ 2 \sum_{ \{u,v\} \in F_3} \phi_{H^*}(\{u,v\}) \cdot X_u X_v. \nonumber  
    \end{align}
        
    Moreover, conditional on $H^*$, the random variables 
    \[
        X_v \sim \Ber(\phi_{H^*}(\{o_{H^*}, v\}))  \qquad \text{for}\,\, v \in \partial B_{H^*}(o_{H^*}, 2),
    \]
     and are mutually independent. Here $\frac{0}{0}$ is interpreted as $0$. 
\end{theorem}

\begin{corollary}\label{cor:cc_local}
    Suppose we are in the set up of Theorem~\ref{thm:cc_local}. Then, the followings hold.
    \begin{enumerate}
        \item [(a)] If $\mu (H^* \,\text{is a tree}) \eqas 1$, then we have the simplification:
        \begin{align*}
            T &= \sum_{\substack{
            u, v \in \cV(\partial B_{H^*}(o_{H^*}, 1))\\
            u \ne v}}
            \phi(d_{H^*}(u), d_{H^*}(v), 1) + 2X + 2 \sum_{\{u,v\} \in F_3}
            \phi (d_{H^*} (u), d_{H^*}(v), 1 ) \cdot X_u X_v.
        \end{align*}
        Here $F_3$ consists of the pairs of vertices of $\partial B_{H^*}(o_{H^*}, 2)$ having a common neighbor in $\partial B_{H^*}(o_{H^*}, 1)$.
        \item[(b)] If $\phi \equiv p$ for some $p \in [0,1]$, then \eqref{eq:cc_local} holds with $T$ replaced by the following quantity
        \begin{align*}
            & T_{H^*}(o_{H^*}) + 2 p |F_1| \\
            & + 
            \frac{2 (|\cE(B_{H^*}(o_{H^*}, 2))| - |\cE(B_{H^*}(o_{H^*}, 1))| - |\cE(\partial B_{H^*}(o_{H^*}, 2))| + p |F_2|) X}{|\cV(\partial B_{H^*}(o_{H^*}, 2))|}  \\
            &+ \frac{2 (|\cE(\partial B_{H^*}(o_{H^*}, 2))| + p |F_3|) X (X-1)}{|\cV(\partial B_{H^*}(o_{H^*}, 2))|(|\cV(\partial B_{H^*}(o_{H^*}, 2))| - 1)},
        \end{align*}
        and conditional on $H^*$,
        \[
            X \sim \Binom (|\cV(\partial B_{H^*}(o_{H^*}, 2))|, p).
        \]
        \item[(c)] If $\mu (H^* \,\text{is a tree}) \eqas 1$ and $\phi \equiv p$ for some $p \in [0,1]$, then
        \begin{equation}\label{eq:cor_local_tree}
            \overline{\cc}_{G_n} \convp \bbE_\mu \left[ \frac{N (N - 1) p + 2 X + P \cdot \frac{X(X-1)}{M(M-1)} p}{N (N - 1)  + (2 N -1)\cdot X + X^2}  \right],
        \end{equation}
        where 
        \begin{equation}\label{eq:cor_cc_local_M_N}
            \begin{split}
            & N:= d_{H^*}(o_{H^*}), \\
            &M:= \sum_{v \in \cV(\partial B_{H^*}(o_{H^*}, 1))} d_{H^*}(v) - N, \\
            & P := \sum_{v \in \cV(\partial B_{H^*}(o_{H^*}, 1))} d_{H^*}(v)^2 - N - 3M,
            \end{split}
        \end{equation}
        and conditional on $H^*$, 
        \begin{align*}
            X \sim \Binom \left(M, p \right).
        \end{align*}
    \end{enumerate}
    Here $\frac{0}{0}$ is interpreted as $0$.
\end{corollary}

\begin{remark}\label{rmk:T_formula}
    In Theorem~\ref{thm:cc_local}, Equation~\eqref{cor:cc_local} holds with $T$ replaced by 
    \[
        \hat{T} := \bbE[T|H^*, X],
    \]
    since the denominator of the RHS is a function of $X$ and not of the individual $X_v$'s. We write an expression of $\hat{T}$ below. To that end, for integers $0 \le s \le t$, let $J_{t,s}$ denote the symmetric polynomial in $t$: $J_{t, s} : \bbR^t \to \bbR$ is defined by
    \begin{equation}\label{eq:J_defn}
        J_{t, s}(x_1, \ldots, x_t) := \sum_{\substack{A \subseteq [t]\\|A| = s}} \prod_{i \in A} x_i \prod_{i \in [t] \setminus A} (1-x_i).
    \end{equation}
    Here $J_{t,0}(x_1, \ldots, x_t)=\prod_{i=1}^t (1-x_i)$, $J_{0,0}=1$, and $J_{t,s}=0$ for $s<0$ or $s>t$. When the number of variables is clear from the context, we just write $J_s$ for $J_{t,s}$.
    Then, we have the following representation of $\hat{T}$:
    \begin{align*}
        \hat{T} &:= T_{H^*}(o_{H^*}) + 2 \sum_{e \in F_1} \phi_{H^*}(e)\\
        & +
        2 \sum_{v \in \cV(\partial B_{H^*}(o_{H^*}, 2))} \left( d_{B_{H^*}(o_{H^*}, 2)} (v) - d_{\partial B_{H^*}(o_{H^*}, 2)} (v) \right) \cdot \bbI(X \ge 1)  \\
        & \qquad \qquad \cdot \frac{\phi_{H^*}(\{o_{H^*}, v\}) J_{X-1}\lf \lf \phi_{H^*}(\{o_{H^*}, u\}) \rf_{u \in \cV( \partial B_{H^*}(o_{H^*}, 2)) \setminus \{v\} } \rf}{J_{X}\lf \lf \phi_{H^*}(\{o_{H^*}, u\}) \rf_{u \in \cV( \partial B_{H^*}(o_{H^*}, 2)) } \rf } \\
        &+ 2 \sum_{\substack{ \{u,v\} \in F_2\\ v \in \cV(\partial B_{H^*}(o_{H^*}, 2))}} \phi_{H^*}(\{u,v\}) \cdot \bbI(X \ge 1)\\
        & \qquad \qquad \cdot\frac{\phi_{H^*}(\{o_{H^*}, v\}) J_{X-1}\lf \lf \phi_{H^*}(\{o_{H^*}, u\}) \rf_{u \in \cV( \partial B_{H^*}(o_{H^*}, 2)) \setminus \{v\} } \rf}{J_{X}\lf \lf \phi_{H^*}(\{o_{H^*}, u\}) \rf_{u \in \cV( \partial B_{H^*}(o_{H^*}, 2)) } \rf } \\
        & + \sum_{\substack{u, v \in \cV(\partial B_{H^*}(o_{H^*}, 2))\\ u \ne v}} \bbI(\{u, v\} \in \cE(H^*)) \cdot \bbI(X \ge 2)\\
        & \qquad \qquad \cdot \frac{\phi_{H^*}(\{o_{H^*}, v\}) \phi_{H^*}(\{o_{H^*}, u\}) J_{X-2}\lf \lf \phi_{H^*}(\{o_{H^*}, w\}) \rf_{w \in \cV( \partial B_{H^*}(o_{H^*}, 2)) \setminus \{u, v\} } \rf}{J_{X}\lf \lf \phi_{H^*}(\{o_{H^*}, w\}) \rf_{w \in \cV( \partial B_{H^*}(o_{H^*}, 2)) } \rf } \\
        &+ 2 \sum_{ \{u,v\} \in F_3} \phi_{H^*}(\{u,v\}) \cdot \bbI(X \ge 2) \\
        & \qquad \qquad \cdot \frac{\phi_{H^*}(\{o_{H^*}, v\}) \phi_{H^*}(\{o_{H^*}, u\}) J_{X-2}\lf \lf \phi_{H^*}(\{o_{H^*}, w\}) \rf_{w \in \cV( \partial B_{H^*}(o_{H^*}, 2)) \setminus \{u, v\} } \rf}{J_{X}\lf \lf \phi_{H^*}(\{o_{H^*}, w\}) \rf_{w \in \cV( \partial B_{H^*}(o_{H^*}, 2)) } \rf }.
    \end{align*}
    One can similarly write an expression for $\hat{T}$ in Corollary~\ref{cor:cc_local} (a). The proof of this easily follows from Theorem~\ref{thm:cc_local} and the following lemma.
\end{remark}

\begin{lemma}\label{lem:conditional_expectation_bernoulli}
    Suppose $\xi_i \sim \Ber(r_i)$ for $i \in [m]$ and are independent. Define $\Xi = \sum_{i=1}^m \xi_i$. Then, we have the following:
    \begin{align}
        &\bbE [\xi_i | \Xi] = \frac{r_i J_{\Xi - 1}\lf \lf r_k \rf_{k \in [m] \setminus \{i\}} \rf}{J_{\Xi}\lf \lf r_k \rf_{k \in [m]} \rf} \bbI (\Xi \ge 1), \label{eq:conditional_xi_1}\\
        &\bbE [\xi_i \xi_j | \Xi] = \frac{r_i r_j J_{\Xi - 2}\lf \lf r_k \rf_{k \in [m] \setminus \{i,j\}} \rf}{J_{\Xi}\lf \lf r_k \rf_{k \in [m]} \rf} \bbI (\Xi \ge 2), \label{eq:conditional_xi_1_xi_2}
    \end{align}
    for $i \in [m]$ and $i, j \in [m]$, $i < j$, respectively, with $r_k \in [0,1]$ for all $k \in [m]$, where for $0 \le s \le t$, $J_{t,s}$ denotes the symmetric polynomial in $t$ variables as defined in \eqref{eq:J_defn}. Here $\frac{0}{0}$ is interpreted as $0$, and the common denominator is positive almost surely. In particular, if $r_i$'s are all equal, then 
    \begin{align}
        &\bbE [\xi_i | \Xi] = \frac{\Xi}{m}, \label{eq:conditional_xi_1_equal}\\
        &\bbE [\xi_i \xi_j | \Xi] = \frac{\Xi (\Xi -1)}{m(m-1)}. \label{eq:conditional_xi_1_xi_2_equal}
    \end{align}
\end{lemma}
The proof of Lemma~\ref{lem:conditional_expectation_bernoulli} is given in Appendix~\ref{sec:aux}.

\begin{corollary}[Degree distribution]\label{cor:degree_distribution}
    Under the assumptions of Theorem~\ref{thm:lwl}, the asymptotic degree distribution of a randomly chosen vertex of $G_n$ is given by the degree distribution of the root of $G^*:=\cT_\phi(H^*)$. More precisely, if conditional on $G_n$, $o_{G_n} \sim \Unif(\cV(G_n))$, then
    \begin{equation}\label{eq:degree_distribution}
        d_{G_n}(o_{G_n}) \convd d_{\overline{H}^*}(o_{\overline{H}^*}) + \overline{X}.
    \end{equation}
    where $\overline{H}^* \sim \overline{\mu}$ and $\overline{\mu}$ is a deterministic measure on $\cG^*$ defined by its action on test functions
    \[
        \bbE_{\overline{\mu}}[h(\overline{H}^*)]  := \bbE [\bbE_\mu [h(H^*)]]  
    \]
    for all $h \in C_b(\cG^*, \rho)$. Moreover conditional on $\overline{H}^*$ with $\overline{H}^* \sim \overline{\mu}$, 
    \[
        \overline{X}:= \sum_{v \in \cV(\partial B_{\overline{H}^*}(o_{\overline{H}^*}, 2))} \overline{X}_v,
    \]
    and
    \[
        \overline{X}_v \sim \Ber(\phi_{\overline{H}^*}(\{o_{\overline{H}^*}, v\})),
        \qquad v \in \cV(\partial B_{\overline{H}^*}(o_{\overline{H}^*}, 2)),
    \]
    are mutually independent. Moreover, for every integer $k \ge 0$,
    \begin{equation}\label{eq:deg_prop}
        \frac{1}{|\cV(G_n)|} \sum_{v \in \cV(G_n)} \bbI(d_{G_n}(v)=k)
        \convp
        \bbE_\mu \left[ \bbP\lt d_{H^*}(o_{H^*}) + X =k \mid H^* \rt \right].
    \end{equation}
\end{corollary}

\begin{theorem}[Global clustering coefficients]\label{thm:cc_global}
    Let $(H_n)_{n \in \bbN}$ be sequence of finite graphs such that $|\cV(H_n)| \to \infty$. Assume that $H_n$ converges locally in probability to $H^* \sim \mu$, where $\mu$ is a possibly random measure on $\cG^*$. Let $\phi$ be a weight function and $G_n:= \cT_\phi(H_n)$ for $n \in \bbN$ and $G^*:= \cT_\phi(H^*)$. Further assume that $\left(|\cV(B_{H_n}(o_{H_n}, 2))|^2 \right)_{n \in \bbN}$ is uniformly integrable, where $o_{H_n} \sim \Unif(\cV(H_n))$, and $\mu(d_{H^*}(o_{H^*}) >1) >0$. Then
    \begin{equation}\label{eq:cc_global}
            \cc_{G_n} \convp \frac{\bbE_\mu [ \overline{T} ]}{\bbE_\mu \left[ \overline{D} \right]},
    \end{equation}
    where
    \begin{align}
        \overline{D} &:= d_{H^*}(o_{H^*}) \cdot (d_{H^*}(o_{H^*}) - 1 ) + 2 d_{H^*}(o_{H^*}) \sum_{v \in \cV(\partial B_{H^*} (o_{H^*}, 2))} \phi_{H^*}(\{o_{H^*}, v\}) \label{eq:D_bar_defn} \\
        &+ \sum_{\substack{
        u, v \in \cV(\partial B_{H^*} (o_{H^*}, 2))\\
        u \ne v}} \phi_{H^*}(\{o_{H^*}, u\}) \cdot \phi_{H^*}(\{o_{H^*}, v\}), \\
        \overline{T} &= T_{H^*}(o_{H^*}) + 2 \sum_{e \in F_1} \phi_{H^*}( e ) \label{eq:T_bar_defn}\\
        &+ 2 \sum_{v \in \cV(\partial B_{H^*}(o_{H^*}, 2))} \left(d_{B_{H^*}(o_{H^*}, 2)} (v) - d_{\partial B_{H^*}(o_{H^*}, 2)} (v) \right) \cdot \phi_{H^*}( \{o_{H^*}, v\} ) \nonumber \\
        &+ 2 \sum_{\substack{\{u,v\} \in F_2 \\ v \in \cV(\partial B_{H^*}(o_{H^*}, 2))}} \phi_{H^*}(\{u,v\}) \cdot \phi_{H^*}(\{o_{H^*}, v\}) \nonumber \\
        &+  \sum_{\substack{u, v \in \cV(\partial B_{H^*}(o_{H^*}, 2))\\ u \ne v}} \bbI(\{u, v\} \in \cE(H^*)) \cdot \phi_{H^*}( \{o_{H^*}, u\} ) \cdot \phi_{H^*}( \{o_{H^*}, v\} ) \nonumber \\
        &+ 2 \sum_{\{u,v\} \in F_3} \phi_{H^*}(\{u,v\}) \cdot \phi_{H^*}(\{o_{H^*}, u\}) \cdot \phi_{H^*}(\{o_{H^*}, v\}). \nonumber
    \end{align}
    Here $F_1$, $F_2$ and $F_3$ are defined before Theorem~\ref{thm:cc_local}.
\end{theorem}

\begin{corollary}\label{cor:cc_global}
    In the set up of Theorem~\ref{thm:cc_global}, the followings hold.
    \begin{enumerate}
        \item[(a)] If $\mu(H^* \, \text{is a tree}) \eqas 1$, then \eqref{eq:cc_global} holds with
        \begin{align*}
            \overline{T} &= \sum_{\substack{u, v \in \cV(\partial B_{H^*}(o_{H^*}, 1))\\ u \ne v}} \phi(d_{H^*}(u), d_{H^*}(v), 1) + 2 \sum_{v \in \partial B_{H^*}(o_{H^*}, 2)} \phi( d_{H^*}(o_{H^*}), d_{H^*}(v), 1 ) \\
            &+ 2 \sum_{\{u,v\} \in F_3} \phi(d_{H^*}(u), d_{H^*}(v), 1)
            \phi (d_{H^*} (o_{H^*}), d_{H^*}(u), 1 )  \phi(d_{H^*}(o_{H^*}), d_{H^*}(v), 1).
        \end{align*}
        \item[(b)] If $\phi \equiv p$ for some $p \in [0,1]$, then \eqref{eq:cc_global} holds with
        \begin{align*}
            \overline{D} &= d_{H^*}(o_{H^*}) \cdot (d_{H^*}(o_{H^*}) - 1 ) + 2 d_{H^*}(o_{H^*}) |\cV(\partial B_{H^*} (o_{H^*}, 2))| p\\
            &+ |\cV(\partial B_{H^*} (o_{H^*}, 2))| (|\cV(\partial B_{H^*} (o_{H^*}, 2))| - 1) p^2,\\
            \overline{T} &= T_{H^*}(o_{H^*}) + 2 p |F_1| + 2 (|\cE(B_{H^*}(o_{H^*}, 2))| - |\cE(B_{H^*}(o_{H^*}, 1))| - |\cE(\partial B_{H^*}(o_{H^*}, 2))|) p \\
            &+ 2 ( |F_2| + |\cE(\partial B_{H^*}(o_{H^*}, 2))|) p^2  + 2 |F_3| p^3.
        \end{align*}
        \item [(c)] If $\mu(H^* \, \text{is a tree}) \eqas 1$ and $\phi \equiv p$ for some $p \in [0,1]$, then, with $M$, $N$, $P$ defined in \eqref{eq:cor_cc_local_M_N},
        \begin{equation}\label{eq:cor_cc_global_tree}
            \cc_{G_n} \convp \frac{\bbE_\mu [N(N-1) + 2M] \cdot p + \bbE_\mu [P] \cdot p^3}{\bbE_\mu \left[N(N-1) \right] + 2\bbE_\mu [M N ] \cdot p + \bbE_\mu [M(M-1)]\cdot p^2}.
        \end{equation}
    \end{enumerate}
\end{corollary}

\begin{remark}
    The proof of Theorem~\ref{thm:cc_global} of course uses Theorem~\ref{thm:lwc_cc_global_remco}. Note that if $\phi >0$, then for the purpose of Theorem~\ref{thm:cc_global}, it is enough to assume that $\mu(|\cV(B_{H^*}(o_{H^*}, 2))| >2) >0$ a.s.. This is because if $|\cV(B_{H^*}(o_{H^*}, 2))| >2$, then with positive $\mu$-probability either $o_{H^*}$ has at least two neighbors or $o_{H^*}$ has at least one vertex at a distance $2$ in $H^*$. In the first case $d_{G^*}(o_{G^*}) \ge 2$ and in the second case, given $H^*$, with positive probability $o_{G^*}$ is joined to this point in $\cV(\partial B_{H^*}(o_{H^*}, 2))$. Thus, obtain that $\bbE_\mu[\bbI(d_{G^*}(o_{G^*}) >1)] >0$ a.s. and Theorem~\ref{thm:lwc_cc_global_remco} applies.
\end{remark}

\begin{remark}\label{rmk:4_moment_UI}
    To ensure that the sequence $\left(|\cV(B_{H_n}(o_{H_n}, 2))|^2 \right)_{n \in \bbN}$ is uniformly integrable, it is enough to assume that the sequence $\left(d_{H_n}(o_{H_n})^4\right)_{n \in \bbN}$ is uniformly integrable. Indeed, for $K\ge1$ and $v\in\cV(H_n)$, define
    \[
        T_{n,K}(v) := \sum_{u \in \cV(\partial B_{H_n}(v, 1))}
            d_{H_n}(u)\bbI_{\{\max\{d_{H_n}(v),d_{H_n}(u)\}>K\}}.
    \]
    Note that
    \[
        |\cV(B_{H_n}(v, 2))|
        \le 1+\sum_{u \in \cV(\partial B_{H_n}(v, 1))}d_{H_n}(u)
        \le 1+K^2+T_{n,K}(v),
    \]
    since the terms not included in $T_{n,K}(v)$ have both endpoint degrees at most $K$. By the Cauchy--Schwarz inequality,
    \begin{align*}
        &\frac{1}{|\cV(H_n)|}\sum_{v\in\cV(H_n)}T_{n,K}(v)^2\\
        &\le \frac{1}{|\cV(H_n)|}\sum_{v\in\cV(H_n)}d_{H_n}(v)
            \sum_{u\in\cV(\partial B_{H_n}(v,1))}
            d_{H_n}(u)^2\bbI_{\{\max\{d_{H_n}(v),d_{H_n}(u)\}>K\}}\\
        &=\frac{2}{|\cV(H_n)|}\sum_{\{u,v\}\in\cE(H_n)}
            d_{H_n}(u)^2d_{H_n}(v) 
            \bbI_{\{\max\{d_{H_n}(v),d_{H_n}(u)\}>K\}}\\
            &\le \frac{2}{3|\cV(H_n)|}\sum_{\{u,v\}\in\cE(H_n)}
            \lf 2d_{H_n}(u)^3 + d_{H_n}(v)^3 \rf
            \bbI_{\{\max\{d_{H_n}(v),d_{H_n}(u)\}>K\}} \\
        &=\frac{1}{|\cV(H_n)|}\sum_{v\in\cV(H_n)}
            d_{H_n}(v)^4\bbI_{\{d_{H_n}(v)>K\}}.
    \end{align*}
    In the penultimate line we have used AM-GM inequality.

    Now, if $|\cV(B_{H_n}(v,2))|>2(1+K^2)$, then $T_{n, K} (v) > 1 + K^2$. The preceding bound then gives
    \[
        |\cV(B_{H_n}(v,2))|<2T_{n,K}(v).
    \]
    Since $o_{H_n}$ is uniformly distributed on $\cV(H_n)$ conditional on $H_n$, we obtain
    \begin{align*}
        &\bbE\lt |\cV(B_{H_n}(o_{H_n},2))|^2
            \bbI_{\{|\cV(B_{H_n}(o_{H_n},2))|>2(1+K^2)\}}\rt\\
        &\le 4\bbE\lt \frac{1}{|\cV(H_n)|}
            \sum_{v\in\cV(H_n)}T_{n,K}(v)^2\rt\\
        &\le 8\bbE\lt \frac{1}{|\cV(H_n)|}
            \sum_{v\in\cV(H_n)}d_{H_n}(v)^4\bbI_{\{d_{H_n}(v)>K\}}\rt\\
        &=8\bbE\lt d_{H_n}(o_{H_n})^4
            \bbI_{\{d_{H_n}(o_{H_n})>K\}}\rt.
    \end{align*}
    Taking the supremum over $n$ and letting $K\to\infty$, the right hand side tends to $0$ by the assumed uniform integrability. This proves uniform integrability of $\left(|\cV(B_{H_n}(o_{H_n},2))|^2\right)_{n\in\bbN}$.
\end{remark}

\section{An empirical study}\label{sec:empirical}

In this section we fit our model to a real world network and see how it competes with its parent. To that end we consider Newman's NetScience coauthorship network, compiled in 2006 from
papers cited in two reviews of network research, with additional references
added by hand \cite{newman2006finding}. In this data each vertex represent a scientist and an edge is put between two vertices if they coauthored at least one paper together. The resulting graph is a simple and undirected. The network has $1589$ many vertices. We do our analysis on its largest connected component which has $n = 379$ many nodes and $m=914$ many edges. The reader can look at \cite{rossi2015network} (\href{https://networkrepository.com/ca-netscience.php}{Network Repository}) for more details about the dataset.

We fit the sparse Erd\H{o}s--R\'enyi model
$\mathrm{ER}(n,\lambda/n)$ and its triadic closure
$\cT_p(\mathrm{ER}(n,\tilde{\lambda}/n))$ with constant closure
probability $p$. For $\mathrm{ER}(n,\lambda/n)$, the maximum-likelihood estimator is
\begin{equation}\label{emp:collaboration_er_mle}
    \hat{\lambda}=\frac{2m}{n-1}=4.8360.
\end{equation}
For the triadic closure model, we use the method of moments, matching mean degree
\[
    \frac{2|E|}{n} = \frac{2m}{n} = 4.8232,
\]
and mean wedge count
\[
    \frac{W}{2n}=\frac{1}{n} \sum_v\binom{d_v}{2}= 16.9314.
\]
Note that in the LHS $2$ appears in the denominator since $W$ as defined in \eqref{eq:wedge_count_defn} represents twice the number of wedges. Using the asymptotic mean degree and mean wedge formulae as derived in subsection~\ref{subsec:er},
the estimating equations are
\begin{align}
 \frac{2 |E|}{n}
 &= \tilde{\lambda}+p\tilde{\lambda}^{\,2},
 \label{emp:collaboration_mean}\\
 \frac{W}{n}
 &= \tilde{\lambda}^{\,2}
    \bigl[1+(\tilde{\lambda}+1)(\tilde{\lambda}p+2)p\bigr].
 \label{emp:collaboration_wedges}
\end{align}
Solving \eqref{emp:collaboration_mean} and \eqref{emp:collaboration_wedges}, we obtain
\[
    \hat{\tilde{\lambda}} =1.8662, \qquad\hat p = 0.8490.
\]
In Table~\ref{emp:collaboration_comparison}, we compare the observed network and the fitted models in terms of a few subgraph count statistics and clustering coefficients. For the fitted models, we use mean of $5000$ simulations to estimate these statistics.  
\begin{table}[htbp]
    \centering
    \small
    \begin{tabular*}{\textwidth}{@{\extracolsep{\fill}}lrrr@{}}
    \toprule
    Statistics & Observed & $\mathrm{ER} \lf n, \frac{\hat{\lambda}}{n} \rf$ & $\cT_{\hat p}\lf\mathrm{ER} \lf n, \frac{\hat{\tilde{\lambda}}}{n}  \rf \rf$\\
    \midrule
    Degree (mean) & 4.823 & 4.824 & 4.769\\
    Degree (variance) & 15.423 & 4.744 & 14.715\\
    Wedges (both open and closed) & 6,417 & 4,399.9 & 6,227.8\\
    Triangles & 921 & 18.7 & 814.1\\
    4-cycles & 2,794 & 67.6 & 2,190.9\\
    4-cliques & 631 & 0.0024 & 326.1\\
    Global clustering coefficient & 0.4306 & 0.0127 & 0.3947\\
    Local clustering coefficient& 0.7412 & 0.0121 & 0.4743\\
    \bottomrule
    \end{tabular*}
    \caption{Comparison of observed and fitted models }
    \label{emp:collaboration_comparison}
\end{table}

As it is evident from Table~\ref{emp:collaboration_comparison}, the triadic closure model substantially improves triangle and other subgraph counts relative to the
Erd\H{o}s-R\'{e}nyi model, although 4-cliques and local clustering remain underestimated. Note that triangle count was not used to estimate the parameters of the model. In order to get better estimates of subgraph counts for denser subgraphs like $4$-cliques, one should consider repeated triadic closure models. We plat to pursue this direction in our future works.

\section{Proofs}\label{sec:proofs}

\subsection{Proof of Theorem~\ref{thm:lwl}}

\begin{proof}[Proof of Theorem~\ref{thm:lwl}]
    Let $G^* := \cT_\phi (H^*)$. It is enough to show that for any integer $r \ge 0$, $K^* \in \cG^*$,
    \begin{equation*}
        p_{G_n}(K^*) \convp \bbE_\mu \lt \bbI(B_{G^*}(o_{G^*}, r ) \cong K^*) \rt.
    \end{equation*}
    Here we write $p_{G_n}(K^*)$ for $p_{G_n}(K^*, r)$. The expectation $\bbE_\mu$ includes the conditional randomness induced by the $Z_{H^*}(e)$'s given $H^*$, with $\mu$ kept fixed when it is random. Since both sides of the above display are in $[0,1]$, this is equivalent to showing that
    \begin{align}
         &\bbE \lt p_{G_n}(K^*) - \bbE_\mu \lt \bbI(B_{G^*}(o_{G^*}, r ) \cong K^*) \rt \rt \to 0, \label{eq:mean_conv}\\
         &\Var \lt p_{G_n}(K^*) - \bbE_\mu \lt \bbI(B_{G^*}(o_{G^*}, r ) \cong K^*) \rt \rt \to 0. \label{eq:var_conv}
    \end{align}
    All sums over $L^*$ below are over isomorphism classes of finite rooted graphs. Now
    \begin{align*}
        p_{G_n}(K^*)
        &= \frac{1}{|\cV(G_n)|} \sum_{v \in \cV(G_n)} \bbI(B_{G_n}(v, r) \cong K^*) \\
        &= \frac{1}{|\cV(H_n)|} \sum_{v \in \cV(H_n)} \sum_{L^* \in \cG^*} \bbI \left(B_{H_n}(v, 2r + 1) \cong L^*, B_{G_n}(v, r) \cong K^* \right).
    \end{align*}
    Fix some labeled rooted graph $L^*$ and a vertex $v \in \cV(H_n)$. Suppose $B_{H_n}(v, 2r + 1) \cong L^*$. Let $\eta_v : \cV(L^*) \to \cV(B_{H_n}(v, 2r + 1))$ be a rooted graph isomorphism.
    Note that each edge of $G_n$ joins vertices at distance at most $2$ in $H_n$. Thus, every vertex of $B_{G_n}(v,r)$ lies in $B_{H_n}(v,2r)$. Moreover, the degrees and common neighbors of vertices of $B_{H_n}(v,2r)$ are completely determined by $B_{H_n}(v,2r+1)$. Hence, the possible added edges with both endpoints in $B_{H_n}(v,2r)$ are determined by $L^*$, and
    \[
        \phi_{H_n}(\eta_v(e)) = \phi_{L^*}(e)
    \]
    for $e=\{w,x\} \in \cW(L^*)$ with $w,x \in \cV(B_{L^*}(o_{L^*},2r))$.
    By the $\phi$-triadic closure operation, a few edges of $\eta_v(\cW(L^*))$ will be added. Denote this set of edges by $\eta_v(E)$ for some $E \subseteq \cW(L^*)$. In order to get $B_{G_n}(v, r) \cong K^*$, $E$ should be such that $B_{L^* \cup E}(o_{L^*}, r) \cong K^*$ where, by an abuse of notation, $L^* \cup E$ denotes the rooted graph obtained by adding the edges $E$ to $L^*$. The other added edges cannot affect this $r$-neighborhood. Thus, on the event $B_{H_n}(v, 2r + 1) \cong L^*$,
    \begin{align*}
        &\{ B_{G_n}(v, r) \cong K^* \} \\
        &= \bigcup_{E \subseteq \cW(L^*)} \{B_{L^* \cup E}(o_{L^*}, r) \cong K^*,\, Z_{\eta_v(e)} =1\,\,\text{for}\,\, e \in E,\, Z_{\eta_v(e)} = 0\,\,\text{for}\,\, e \in \cW(L^*) \setminus E \}.
    \end{align*}
    Here $Z_{\eta_v(e)}$ denotes $Z_{H_n}(\eta_v(e))$. Then, we have
    \begin{align*}
        p_{G_n}(K^*) &= \frac{1}{|\cV(H_n)|} \sum_{v \in \cV(H_n)} \sum_{L^* \in \cG^*} \bbI \left(B_{H_n}(v, 2r + 1) \cong L^* \rf \\
        & \cdot \sum_{E \subseteq \cW(L^*)} \bbI \lf B_{L^* \cup E}(o_{L^*}, r) \cong K^* \rf \cdot \prod_{e \in E} Z_{H_n}(\eta_v(e)) \cdot \prod _{e \in \cW(L^*) \setminus E} (1 - Z_{H_n} (\eta_v(e))).
    \end{align*}
    Now, for $L^* \in \cG^*$ and $E \subseteq \cW(L^*)$, define
    \begin{align*}
        & T_{H_n, v}(L^*, E) := \prod_{e \in E} Z_{H_n} (\eta_v(e)) \cdot \prod _{e \in \cW(L^*) \setminus E} (1 - Z_{H_n} (\eta_v(e))),\\
        & R_{H_n, v}(L^*) := \sum_{E \subseteq \cW(L^*)} \bbI \lf B_{L^* \cup E}(o_{L^*}, r) \cong K^* \rf \cdot T_{H_n, v}(L^*, E), \\
        & \overline{T}(L^*, E) := \prod_{e \in E} \phi_{L^*}(e)  \cdot \prod _{e \in \cW(L^*) \setminus E} (1 - \phi_{L^*}(e)), \\
        & \overline{R}(L^*) := \sum_{E \subseteq \cW(L^*)} \bbI \lf B_{L^* \cup E}(o_{L^*}, r) \cong K^* \rf \cdot \overline{T}(L^*, E).
    \end{align*}
    We then have
    \[
        p_{G_n}(K^*) = \frac{1}{|\cV(H_n)|} \sum_{v \in \cV(H_n)} \sum_{L^* \in \cG^*} \bbI \left(B_{H_n}(v, 2r + 1) \cong L^* \rf \cdot R_{H_n, v}(L^*).
    \]
    Note that $T_{H_n,v}(L^*,E)$ and $R_{H_n,v}(L^*)$ are well defined on the event $B_{H_n}(v, 2r + 1) \cong L^*$, and we set them equal to $0$ otherwise. The indicator $\bbI(B_{L^* \cup E}(o_{L^*},r) \cong K^*)$ depends only on the edges of $E$ with both endpoints in $B_{L^*}(o_{L^*},2r)$. Summing over the choices of all other edges in $R_{H_n,v}(L^*)$ and $\overline{R}(L^*)$ gives factors $Z_{H_n} (\eta_v(e))+(1-Z_{H_n} (\eta_v(e)))=1$ and $\phi_{L^*}(e)+(1-\phi_{L^*}(e))=1$, respectively. For the remaining edges the closure probabilities agree, as observed above. Moreover, the indicator of the event $B_{H_n}(v, 2r + 1) \cong L^*$ is $H_n$ measurable and hence
    \[
        \bbE \lt \bbI \left(B_{H_n}(v, 2r + 1) \cong L^* \rf \cdot R_{H_n, v}(L^*) \mid H_n \rt = \bbI \left(B_{H_n}(v, 2r + 1) \cong L^* \rf \cdot \overline{R}(L^*).
    \]
    Define $f : \cG^* \to [0,1]$ by
    \[
        f(M^*) := \sum_{L^* \in \cG^*} \overline{R}(L^*) \cdot \bbI \lf B_{M^*}(o_{M^*}, 2r+1) \cong L^* \rf.
    \]
    Clearly, $f$ is continuous, since it depends only on the $(2r+1)$-neighborhood of the root. Similarly,
    \begin{align*}
        \bbE_\mu \lt \bbI(B_{G^*}(o_{G^*}, r ) \cong K^*) \rt
        &= \sum_{L^* \in \cG^*} \mu \lf B_{H^*}(o_{H^*}, 2r + 1) \cong L^* \rf \cdot \overline{R}(L^*)\\
        &= \bbE_\mu \lt f(H^*) \rt.
    \end{align*}

    \textbf{Proof of \eqref{eq:mean_conv}.} By the tower property of conditional expectation,
    \begin{align*}
        &\bbE \lt p_{G_n}(K^*) - \bbE_\mu \lt f(H^*) \rt \rt\\
        &= \bbE \lt \bbE \lt p_{G_n}(K^*) \mid H_n \rt - \bbE_\mu \lt f(H^*) \rt \rt\\
        &= \bbE \lt \frac{1}{|\cV(H_n)|} \sum_{v \in \cV(H_n)} f((H_n, v )) - \bbE_\mu \lt f(H^*) \rt \rt.
    \end{align*}
    Since $H_n$ converges locally in probability to $H^*$,
    \[
        \frac{1}{|\cV(H_n)|} \sum_{v \in \cV(H_n)} f((H_n, v )) \convp \bbE_\mu \lt f(H^*) \rt.
    \]
    Since both sides are bounded by $1$, the difference converges to $0$ in $L^1$ and in $L^2$. This proves equation~\eqref{eq:mean_conv}.

    \textbf{Proof of \eqref{eq:var_conv}.} Now using the conditional variance formula and the variance bound for a sum,
    \begin{align*}
        &\Var \lt p_{G_n}(K^*) - \bbE_\mu \lt f(H^*) \rt \rt\\
        &\le 2\Var \lt \bbE \lt p_{G_n}(K^*) \mid H_n \rt - \bbE_\mu \lt f(H^*) \rt \rt
        + 2\Var \lt p_{G_n}(K^*) - \bbE \lt p_{G_n}(K^*) \mid H_n \rt \rt\\
        &= 2\Var \lt \frac{1}{|\cV(H_n)|} \sum_{v \in \cV(H_n)} f((H_n, v )) - \bbE_\mu \lt f(H^*) \rt \rt
        + 2\bbE \lt \Var \lt p_{G_n}(K^*) \mid H_n \rt \rt.
    \end{align*}
    We show that each term in the above display converges to $0$. For the first term, the $L^2$ convergence proved above gives
    \[
        \Var \lt \frac{1}{|\cV(H_n)|} \sum_{v \in \cV(H_n)} f((H_n, v )) - \bbE_\mu \lt f(H^*) \rt \rt \to 0.
    \]
    Next we consider the other part:
    \begin{align*}
        &\bbE \lt \Var \lt p_{G_n}(K^*) \mid H_n \rt \rt\\
        &= \bbE \lt \Var \lt \frac{1}{|\cV(H_n)|} \sum_{v \in \cV(H_n)} \sum_{L^* \in \cG^*} \bbI \left(B_{H_n}(v, 2r + 1) \cong L^* \rf \cdot R_{H_n, v}(L^*) \bigg| H_n \rt \rt.
    \end{align*}
    Define for $u, v \in \cV(H_n)$,
    \begin{align*}
        &P_{H_n}(v) := \Var \lt \sum_{L^* \in \cG^*} \bbI \left(B_{H_n}(v, 2r + 1) \cong L^* \rf \cdot R_{H_n, v}(L^*) \bigg| H_n \rt, \\
        &Q_{H_n}(u, v) := \Cov \bigg[ \sum_{L^* \in \cG^*} \bbI \left(B_{H_n}(u, 2r + 1) \cong L^* \rf \cdot R_{H_n, u}(L^*),\\
        & \qquad \qquad \qquad \sum_{L^* \in \cG^*} \bbI \lf B_{H_n}(v, 2r + 1) \cong L^* \rf \cdot R_{H_n, v}(L^*) \bigg| H_n \bigg].
    \end{align*}
    Clearly, $\left| P_{H_n}(v) \right|, \left| Q_{H_n} (u,v) \right| \le 1$. Then,
    \begin{align*}
        &\bbE \lt \Var \lt p_{G_n}(K^*) \mid H_n \rt \rt \\
        &= \bbE \lt \frac{1}{|\cV(H_n)|^2} \sum_{v \in \cV(H_n)} P_{H_n}(v)
        + \frac{1}{|\cV(H_n)|^2} \sum_{\substack{u, v \in \cV(H_n)\\ u \ne v}} Q_{H_n}(u, v) \rt.
    \end{align*}
    The first term in the above display is bounded by $\bbE \lt \frac{1}{|\cV(H_n)|} \rt$, which converges to $0$. For the second term, note that if $\dist_{H_n}(u,v) \ge 4r+2$, then $B_{H_n}(u,2r)$ and $B_{H_n}(v,2r)$ have disjoint vertex sets. Thus, the corresponding indicators depend on disjoint collections of $Z_{H_n}(e)$'s, which are independent conditional on $H_n$, and $Q_{H_n}(u,v) = 0$. Hence, the second term, in absolute value, is bounded by
    \[
        \bbE \lt \frac{|\{(u, v) \in \cV(H_n)^2: \dist_{H_n}(u, v) \le 4r+1\}|}{|\cV(H_n)|^2} \rt.
    \]
    The quantity inside the expectation is precisely the proportion of vertex pairs in $H_n$ for which the distance is less than or equal to $4r+1$. By Corollary~2.20 of \cite{vanderhofstad2024random}, this converges to $0$ in probability. Since the proportion is also bounded by $1$, its expectation converges to $0$. This finally yields
    \[
        \bbE \lt \Var \lt p_{G_n}(K^*) \mid H_n \rt \rt \to 0,
    \]
    and proves \eqref{eq:var_conv}. Together with \eqref{eq:mean_conv}, this finishes the proof of Theorem~\ref{thm:lwl}.
\end{proof}

\subsection{Proof of Theorems~\ref{thm:cc_local} and \ref{thm:cc_global}}

\begin{proof}[Proof of Theorem~\ref{thm:cc_local}]   
    By Theorems~\ref{thm:lwl} and \ref{thm:lwc_cc_local_remco},
    \begin{align*}
        \overline{\cc}_{G_n} \convp \bbE_\mu \left[ \frac{T_{G^*}(o_{G^*}) }{d_{G^*}(o_{G^*})(d_{G^*}(o_{G^*}) - 1)} \right].
    \end{align*}
    We now compute distribution of the quantity inside the expectation. Note that $B_{H^*}(o_{H^*},1) \subseteq B_{G^*}(o_{G^*}, 1)$. Moreover,
    \[
        d_{G^*}(o_{G^*}) = d_{H^*}(o_{H^*}) + X.
    \]
    So, 
    \begin{align*}
        d_{G^*}(o_{G^*})(d_{G^*}(o_{G^*}) - 1)
        &= d_{G^*}(o_{G^*})^2 - d_{G^*}(o_{G^*}) \\
        &= d_{H^*}(o_{H^*}) \cdot (d_{H^*}(o_{H^*}) - 1)  + (2 d_{H^*}(o_{H^*}) -1)\cdot X + X^2. 
    \end{align*}
    We now compute $T_{G^*}(o_{G^*})$. As observed, there are $d_{H^*}(o_{H^*}) + X$ many neighbors of $o_{G^*}$ in $G^*$. Triangles are formed by joining these vertices among themselves. This can happen in three ways. We compute these counts separately.

    Firstly, two vertices of $\partial B_{H^*}(o_{H^*}, 1)$ can be joined to make a triangle. Call these the type-$1$ triangles. There are $\frac{1}{2} T_{H^*}(o_{H^*})$ such triangles inherited from $H^*$. The contribution in $T_{G^*}(o_{G^*})$ from the type-$1$ triangles is precisely 
    \[
        T^{(1)} := T_{H^*}(o_{H^*}) + 2 \sum_{e \in F_1} Z_{H^*}(e) .
    \]

    Secondly, consider the triangles whose one endpoint lies in $\partial B_{H^*}(o_{H^*}, 1)$ and the other in $\partial B_{H^*}(o_{H^*}, 2)$. Call these the type-$2$ triangles. Note that for each $v \in \partial B_{H^*}(o_{H^*}, 2)$, if $o_{H^*}$ and $v$ are connected in the triadic closure, each edge in $G^*$ connecting $v$ to a vertex of $\partial B_{H^*}(o_{H^*}, 1)$ gives rise to a triangle. There are $\left( d_{B_{H^*}(o_{H^*}, 2)} (v) - d_{\partial B_{H^*}(o_{H^*}, 2)} (v) \right)$ many such edges inherited from $H^*$. Moreover, for $u \in \partial B_{H^*}(o_{H^*}, 1)$ with $\{u,v\} \in \cW(H^*)$, the edge $\{u,v\}$ is added in the triadic closure if $Z_{H^*}(\{u,v\})=1$. Thus, the contribution in $T_{G^*}(o_{G^*})$ from the type-$2$ triangles is precisely
    \begin{align*}
        T^{(2)} &:= 2 \sum_{v \in \partial B_{H^*}(o_{H^*}, 2)} \left( d_{B_{H^*}(o_{H^*}, 2)} (v) - d_{\partial B_{H^*}(o_{H^*}, 2)} (v) \right) \cdot X_v \\
        & \qquad + 2 \sum_{\substack{\{u,v\} \in F_2 \\ v \in \partial B_{H^*}(o_{H^*}, 2)}} Z_{H^*}(\{u,v\}) \cdot X_v.
    \end{align*}

    Lastly, the type-$3$ triangles are those whose both endpoints come from $\partial B_{H^*}(o_{H^*}, 2)$. Note that for distinct $u, v \in \partial B_{H^*}(o_{H^*}, 2)$, a triangle is formed in $G^*$ with vertices $o_{H^*}, u, v$ if $X_u=X_v=1$ and $\{u,v\} \in \cE(G^*)$. The edge $\{u,v\}$ can either be inherited from $H^*$ or added in the triadic closure. In the latter case, $\{u,v\} \in \cW(H^*)$ and $Z_{H^*}(\{u,v\})=1$. Hence, the contribution in $T_{G^*}(o_{G^*})$ from the triangles of type-$3$ is given by
    \begin{align*}
        T^{(3)} &:= \sum_{\substack{u, v \in \partial B_{H^*}(o_{H^*}, 2)\\ u \ne v}} \bbI(\{u, v\} \in \cE(H^*)) \cdot X_u X_v \\
        & \qquad + 2 \sum_{\{u,v\} \in F_3} Z_{H^*}(\{u,v\}) \cdot X_u X_v.
    \end{align*}
    Hence
    \begin{equation*}
        \begin{split}
            T_{G^*}(o_{G^*}) &= T^{(1)} + T^{(2)} + T^{(3)},
        \end{split}
    \end{equation*}
    and consequently,
    \[
        \overline{\cc}_{G_n} \convp \bbE_\mu \left[ \frac{T_{G^*}(o_{G^*})}{D} \right].
    \]
    Note that, conditional on $H^*$, $Z_{H^*}(e)$ variables are independent. Thus, we have
    \begin{align*}
        &\bbE \lt Z_{H^*}(e) | (X_v)_{v \in \cV(\partial B_{H^*}(o_{H^*}, 2))}, H^* \rt = \phi_{H^*} (e)\qquad \text{for}\,\, e \in F_1 \cup F_2 \cup F_3.
    \end{align*}
    Moreover conditional on $H^*$, $D$ is $(X_v)_{v \in \cV(\partial B_{H^*}(o_{H^*}, 2))}$ measurable. So, by the tower property of conditional expectation, we can write
    \begin{align*}
        \bbE_\mu \left[ \frac{T_{G^*}(o_{G^*})}{D} \right] &= \bbE_\mu  \lt \frac{\bbE_{\mu} \left[ T_{G^*}(o_{G^*}) | (X_v)_{v \in \cV(\partial B_{H^*}(o_{H^*}, 2))}, H^* \right]}{D} \rt \\
        &= \bbE_\mu \lt \frac{T}{D} \rt,
    \end{align*}
    This finishes the proof of Theorem~\ref{thm:cc_local}.
\end{proof}

\begin{proof}[Proof of Corollary~\ref{cor:cc_local}]
    In order to prove (a), first note that when $H^*$ is a tree, $T_{H^*}(o_{H^*}) =0$. Further, since there is no cycle in $H^*$, for $v \in \partial B_{H^*}(o_{H^*}, 2)$, $d_{B_{H^*}(o_{H^*}, 2)} (v) =1$, $d_{\partial B_{H^*}(o_{H^*}, 2)} (v) =0$ and there are no edges of $H^*$ between two vertices of $\partial B_{H^*}(o_{H^*}, 2)$. Also, $F_2=\emptyset$, since the length of any path between a vertex of $\partial B_{H^*}(o_{H^*}, 1)$ and a vertex of $\partial B_{H^*}(o_{H^*}, 2)$ is either $1$ or $3$. Two distinct vertices of $\partial B_{H^*}(o_{H^*}, 2)$ are at distance $2$ precisely when they have a common neighbor in $\partial B_{H^*}(o_{H^*}, 1)$. This completes the proof of (a).

    To prove (b), one can consider each term in the expression of $T$, find its expectation conditional on $X$ and simplify it using Lemma~\ref{lem:conditional_expectation_bernoulli}. 
    
    Lastly, to prove (c), we use (b). Moreover since $H^*$ is a tree,
    \begin{align*}
        & T_{H^*}(o_{H^*}) = 0, \\
        &|F_1|=\frac{N(N-1)}{2}, \\
        & |\cV(\partial B_{H^*}(o_{H^*}, 2))| = M,\\
        &|\cE(B_{H^*}(o_{H^*}, 2))|  = M+ N, \\ &|\cE(B_{H^*}(o_{H^*}, 1))| = N ,\\
        & |\cE(\partial B_{H^*}(o_{H^*}, 2))| = 0, \\
        &|F_2| = 0, \\
        & |F_3| = \sum_{v \in \cV(\partial B_{H^*}(o_{H^*}, 1))} \binom{d_{H^*}(v)-1}{2} = \frac{P}{2},
    \end{align*}
    and hence the proof of (c) is complete.
\end{proof}

\begin{proof}[Proof of Corollary~\ref{cor:degree_distribution}]
    By Theorem~\ref{thm:lwl}, $G_n$ converges locally in probability to $G^*:=\cT_\phi(H^*)$. For every integer $k \ge 0$, the function
    \[
        f_k(K^*):=\bbI(d_{K^*}(o_{K^*})=k), \qquad K^* \in \cG^*,
    \]
    is bounded and continuous, since it depends only on $B_{K^*}(o_{K^*},1)$. Hence, by Definition~\ref{def:lwc}(b),
    \begin{align*}
        \frac{1}{|\cV(G_n)|} \sum_{v \in \cV(G_n)} \bbI(d_{G_n}(v)=k)
        &\convp
        \bbE_\mu \left[ \bbP [ d_{G^*}(o_{G^*})=k \mid H^* ] \right] \\
        &= \bbE_\mu \left[ \bbP [d_{H^*}(o_{H^*}) + X = k \mid H^* ] \right].
    \end{align*}
    This proves \eqref{eq:deg_prop}. Both sides of \eqref{eq:deg_prop} lie in $[0,1]$, so convergence in probability also gives convergence in $L^1$. Taking expectations, we obtain
    \begin{align*}
        \bbP(d_{G_n}(o_{G_n})=k)
        &= \bbE \left[ \frac{1}{|\cV(G_n)|} \sum_{v \in \cV(G_n)} \bbI(d_{G_n}(v)=k) \right] \\
        &\to \bbE \left[ \bbE_\mu \left[ \bbP [ d_{H^*}(o_{H^*}) + X=k \mid H^* ] \right] \right] \\
        &= \bbP_{\overline{\mu}}(d_{\overline{H}^*}(o_{\overline{H}^*})+ \overline{X}=k).
    \end{align*}
    Since $H^*$ is locally finite, so is $\overline{H}^*$ and hence $d_{\overline{H}^*}(o_{\overline{H}^*})+\overline{X}$ is finite almost surely. Therefore, convergence of these probabilities for every integer $k \ge 0$ proves \eqref{eq:degree_distribution}.
\end{proof}

\begin{proof}[Proof of Theorem~\ref{thm:cc_global}]
    Take $o_{G_n}=o_{H_n}$. Note that
    \[
        \bbE_\mu [ \bbI(d_{G^*}(o_{G^*}) > 1)] \ge \mu (d_{H^*}(o_{H^*})>1) >0,
    \]
    almost surely. Moreover, since $\cV(B_{G_n}(v,1)) \subseteq \cV(B_{H_n}(v,2))$ for all $v \in \cV(H_n)$,
    \[
        d_{G_n}(o_{G_n})^2 \le \left(|\cV(B_{H_n}(o_{H_n}, 2))| - 1\right)^2 \le |\cV(B_{H_n}(o_{H_n}, 2))|^2,
    \]
    and hence $\left(d_{G_n}(o_{G_n})^2\right)_{n \in \bbN}$ is uniformly integrable. By Theorem~\ref{thm:lwl}, $G_n$ converges locally in probability to $G^*$. So, an use of Theorem~\ref{thm:lwc_cc_global_remco} yields
    \begin{equation*}
            \cc_{G_n} \convp \frac{\bbE_\mu [ T_{G^*}(o_{G^*}) ]}{\bbE_\mu[d_{G^*}(o_{G^*})(d_{G^*}(o_{G^*}) - 1)]}.
    \end{equation*}
    Now, \eqref{eq:cc_global} follows by noting that $\bbE_\mu [T_{G^*}(o_{G^*})] = \bbE_\mu [T] = \bbE_\mu [\bbE_\mu [T|H^*]]$ and $\bbE_\mu [T | H^*] = \overline{T}$ where $T$ is defined in \eqref{eq:T_defn}, and similarly for the denominator.  
\end{proof}

\begin{proof}[Proof of Corollary~\ref{cor:cc_global}]
    The proof essentially follows from Theorem~\ref{thm:cc_global} along with the counting arguments in the proof Corollary~\ref{cor:cc_local}. In order to prove (a), we use the tree properties noted there to get $T_{H^*}(o_{H^*}) =0$, $F_2=\emptyset$ and $\cE(\partial B_{H^*}(o_{H^*}, 2))=\emptyset$. Moreover, $d_{B_{H^*}(o_{H^*}, 2)}(v)-d_{\partial B_{H^*}(o_{H^*}, 2)}(v)=1$ for $v \in \partial B_{H^*}(o_{H^*}, 2)$, and $\phi_{H^*}(\{u,v\})=\phi(d_{H^*}(u),d_{H^*}(v),1)$ for $\{u,v\} \in \cW(H^*)$. Substituting these identities in \eqref{eq:T_bar_defn} gives the stated form of $\overline{T}$. This completes the proof of (a).

    To prove (b), one can consider each term in the expressions of $\overline{D}$ and $\overline{T}$ and simplify it using $\phi \equiv p$. Note that
    \begin{align*}
        &2 |F_1| = d_{H^*}(o_{H^*}) \cdot (d_{H^*}(o_{H^*}) - 1) - T_{H^*}(o_{H^*}),\\
        &\sum_{v \in \partial B_{H^*}(o_{H^*}, 2)} \left(d_{B_{H^*}(o_{H^*}, 2)} (v) - d_{\partial B_{H^*}(o_{H^*}, 2)} (v) \right)\\
        &\qquad = |\cE(B_{H^*}(o_{H^*}, 2))| - |\cE(B_{H^*}(o_{H^*}, 1))| - |\cE(\partial B_{H^*}(o_{H^*}, 2))|,\\
        &\sum_{\substack{u, v \in \partial B_{H^*}(o_{H^*}, 2)\\ u \ne v}} \bbI(\{u,v\} \in \cE(H^*)) = 2 |\cE(\partial B_{H^*}(o_{H^*}, 2))|.
    \end{align*}
    Each pair in $F_2$ contributes $2p^2$ to $\overline{T}$ and each pair in $F_3$ contributes $2p^3$ to $\overline{T}$. This completes the proof of (b).

    Lastly, to prove (c), we use (b) along with the counting identities in the proof of the corollary to Theorem~\ref{thm:cc_local}. In particular, $T_{H^*}(o_{H^*}) = 0$, $|F_1|=N(N-1)/2$, $|\cV(\partial B_{H^*}(o_{H^*}, 2))| = M$, $|F_2| = 0$ and $|F_3| = P/2$.
    Consequently,
    \begin{align*}
        \overline{D} &= N(N-1) + 2MNp + M(M-1)p^2,\\
        \overline{T} &= (N(N-1) + 2M)p + Pp^3,
    \end{align*}
    and hence the proof of (c) is complete.
\end{proof}

\section{Examples}\label{sec:examples}

In this section, we consider few classical examples including sparse Erd\H{o}s-R\'{e}nyi random graphs, random $d$-regular graphs, configuration models, preferential attachment models. Denote them by $H_n$ and their triadic closures by $G_n = \cT_p(H_n)$ for some $p\in [0,1]$. We describe the local limit of these graphs, denoted by $H^*$. Then, Theorem~\ref{thm:lwl} directly applies to these models to yield the local limit of the triadic closures $G^* = \cT_p(H^*)$. In all these examples, the local limit of the parent graph is a tree, $\overline{\cc}_{H^*} =0$ and $\cc_{H^*} =0$. However, by its very construction, the triadic closures contain triangles with positive probability provided $p>0$. By Corollary~\ref{cor:cc_local}(c) and Corollary~\ref{cor:cc_global}(c), we get that 
\begin{align*}
    & \overline{\cc}_{G_n} \convp \overline{\cc}_{G^*},  & \cc_{G_n} \convp \cc_{G^*},
\end{align*}
respectively. In order to describe $\overline{\cc}_{G^*}$, we need the joint distribution of $M$, $N$ and $P$. To compute $\cc_{G^*}$, we evaluate $\bbE[N]$, $\bbE[N^2]$, $\bbE[M]$, $\bbE[MN]$, $\bbE[M^2]$, $\bbE[P]$. Moreover, by Remark~\ref{rmk:4_moment_UI}, it is sufficient to verify that $\sup_n \bbE \lt d_{H_n}(o_{H_n})^{4+\epsilon} \rt <\infty$ for some $\epsilon>0$. Further, Corollary~\ref{cor:degree_distribution} provides the asymptotic degree distribution of $G_n$. We work out the details below.

\subsection{Erd\H{o}s-R\'{e}nyi random graphs}\label{subsec:er}

Recall that \emph{Erd\H{o}s-R\'{e}nyi random graph} with $n$ vertices and connection probability $q$, denoted by $\ER(n, q)$, is a random graph formed by taking $n$ vertices and joining each of the $\binom{n}{2}$ many edges independently with probability $q$ where $q \in [0,1]$. We shall consider sparse Erd\'{o}s-R\'{e}nyi model. More precisely, assume that $H_n \sim \ER \left( n, \frac{\lambda}{n} \right)$, where $\lambda >0$ is a constant representing the asymptotic average degree in the graph. It is well known (e.g. look at Theorem~2.18 of \cite{vanderhofstad2024random}) that $H_n$ converges locally in probability to the \emph{Poisson Galton-Watson tree} with parameter $\lambda$. We denote it by $\PGWT(\lambda)$. So, $H^* \sim \PGWT(\lambda)$. So, $G^* = \cT_p(\PGWT(\lambda))$.

By definition of $\PGWT$, $N \sim \Poisson(\lambda)$ and conditional on $N$, $M = \sum_{i=1}^N (d_{H^*}(v_i)-1)$ where $v_1, \ldots, v_N$ are the neighbors of $o_{H^*}$ and $d_{H^*}(v_i)-1 \overset{\iid}{\sim} \Poisson(\lambda)$ for $i \in [N]$. So, $M | N \sim \Poisson(N \lambda)$. Moreover, by \eqref{eq:cor_cc_local_M_N},
\[
    P = \sum_{i=1}^N d_{H^*}(v_i)^2 - N - 3M = \sum_{i=1}^N (d_{H^*}(v_i)-1)(d_{H^*}(v_i)-2).
\]
This gives the joint distribution of $M$, $N$ and $P$ and lets us calculate $\overline{\cc}_{G^*}$.

\begin{figure}[htbp]
    \centering
    \includegraphics[width=0.6\linewidth]{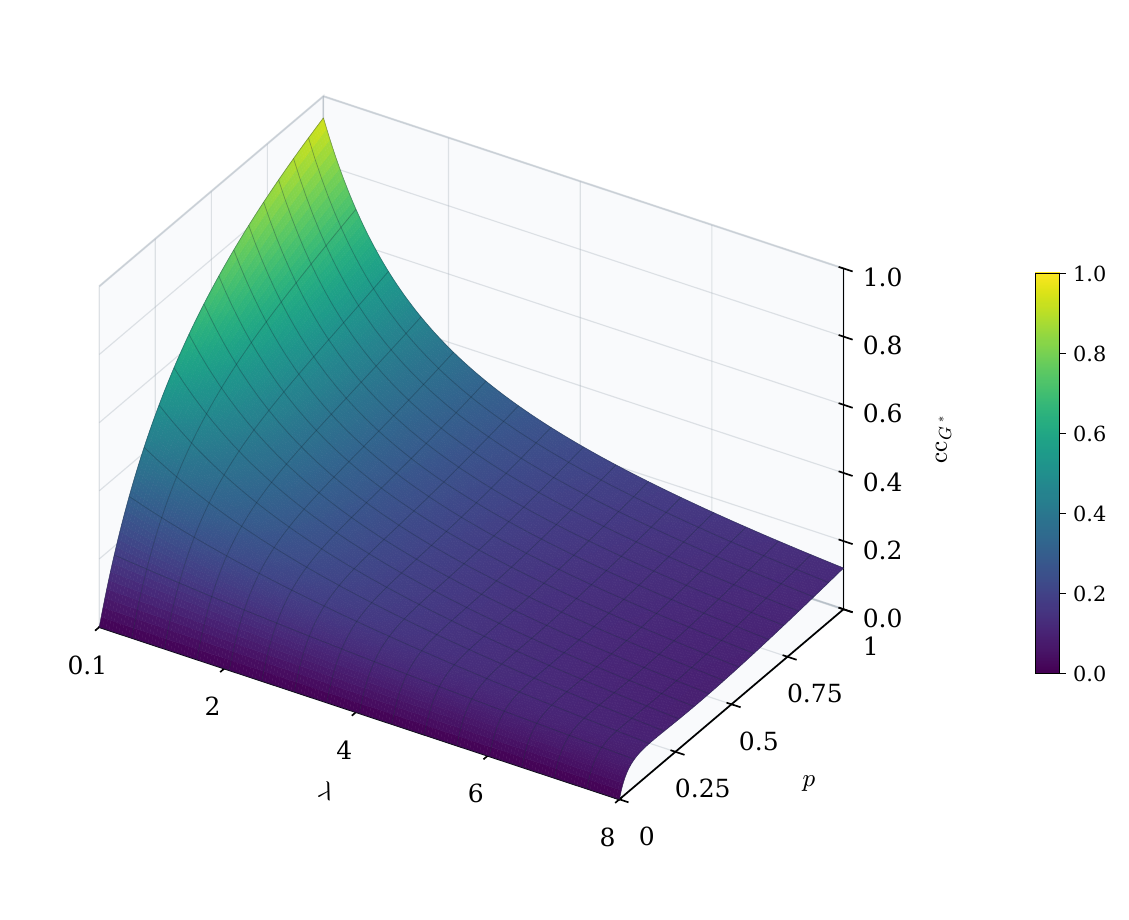}
    \caption{Limiting global clustering coefficient of
    $\cT_p(\mathrm{ER}(n,\lambda/n))$.}
    \label{fig:cc_global_er_surface}
\end{figure}

In order to apply Corollary~\ref{cor:cc_global}, first note that
\begin{align*}
    \bbE \lt d_{H_n}(o_{H_n})^5 \rt &= \bbE \lt d_{H_n}(v)^5 \rt,
\end{align*}
for any fixed $v \in \cV(H_n)$. Moreover, $d_{H_n}(v) \sim \Binom(n-1, \frac{\lambda}{n})$. From moment formula of Binomial distribution, one obtains
\[
    \lim_n \bbE \lt d_{H_n}(v)^5 \rt = \lambda^5 + 10 \lambda^4+ 25 \lambda^3 + 15 \lambda^2 + \lambda.
\] 
We have thus verified that $\sup_n \bbE \lt d_{H_n}(o_{H_n})^5 \rt < \infty$. Now, one can now compute that 
\begin{align*}
    &\bbE[N] = \lambda,\\ 
    &\bbE[N^2] = \lambda^2 + \lambda,\\
    &\bbE[M] = \bbE[\bbE[M|N]] = \bbE[N \lambda] = \lambda^2,\\ 
    &\bbE[MN] = \bbE[N \bbE[M|N]] = \bbE[N \lambda\cdot N] = \lambda(\lambda^2 + \lambda) = \lambda^2(\lambda + 1),\\
    &\bbE[M^2] = \bbE[(N \lambda)^2 + (N \lambda)] = (\lambda^2 + \lambda) \cdot \lambda^2 + \lambda \cdot \lambda = \lambda^2(\lambda^2 + \lambda +1),\\
    &\bbE[P] = \bbE[\bbE[P|N]] = \bbE[N \lambda^2] = \lambda^3.
\end{align*}
Thus,
\begin{equation}\label{eq:cc_global_er}
    \begin{split}
    \cc_{G^*} &= \frac{((\lambda^2 + \lambda) - \lambda + 2 \lambda^2)\cdot p + \lambda^3 \cdot p^3}{(\lambda^2 + \lambda) - \lambda + 2 \lambda^2 (\lambda + 1)\cdot p + (\lambda^2(\lambda^2 + \lambda +1) - \lambda^2)\cdot p^2}\\
    &= \frac{3 p + \lambda p^3}{1 + (\lambda + 1) (\lambda p + 2) p}.
    \end{split}
\end{equation}
In Figure~\ref{fig:cc_global_er_surface}, we plot $\cc_{G^*}$ as a function of $\lambda$ and $p$.

\subsection{Random $k$-regular graphs} 

Fix $k \in \bbN$ with $k \ge 2$. Suppose for $n>k$ and $nk$ even, $H_n$ is chosen uniformly at random from the set of $k$-regular graphs on $n$ vertices. It is well known (e.g. see Theorem~2.17 of \cite{vanderhofstad2024random}) that $H_n$ converges locally in probability to the rooted $k$-regular tree which we denote by $H^*$. In this case, $M \equiv k(k-1)$ and $N \equiv k$. Moreover,
\[
    P = \sum_{v \in \cV(\partial B_{H^*}(o_{H^*}, 1))} (d_{H^*}(v)-1)(d_{H^*}(v)-2) \equiv k(k-1)(k-2).
\]
So, by Corollary~\ref{cor:cc_local}(c),
\[
    \overline{\cc}_{G^*} = \bbE \lt \frac{k(k-1)p + 2X + \frac{(k-2)p X(X-1)}{k(k-1)-1}}{k(k-1) + (2k-1)X + X^2} \rt,
\]
where $X \sim \Binom(k(k-1), p)$.

\begin{figure}[htbp]
    \centering
    \includegraphics[width=0.6\linewidth]{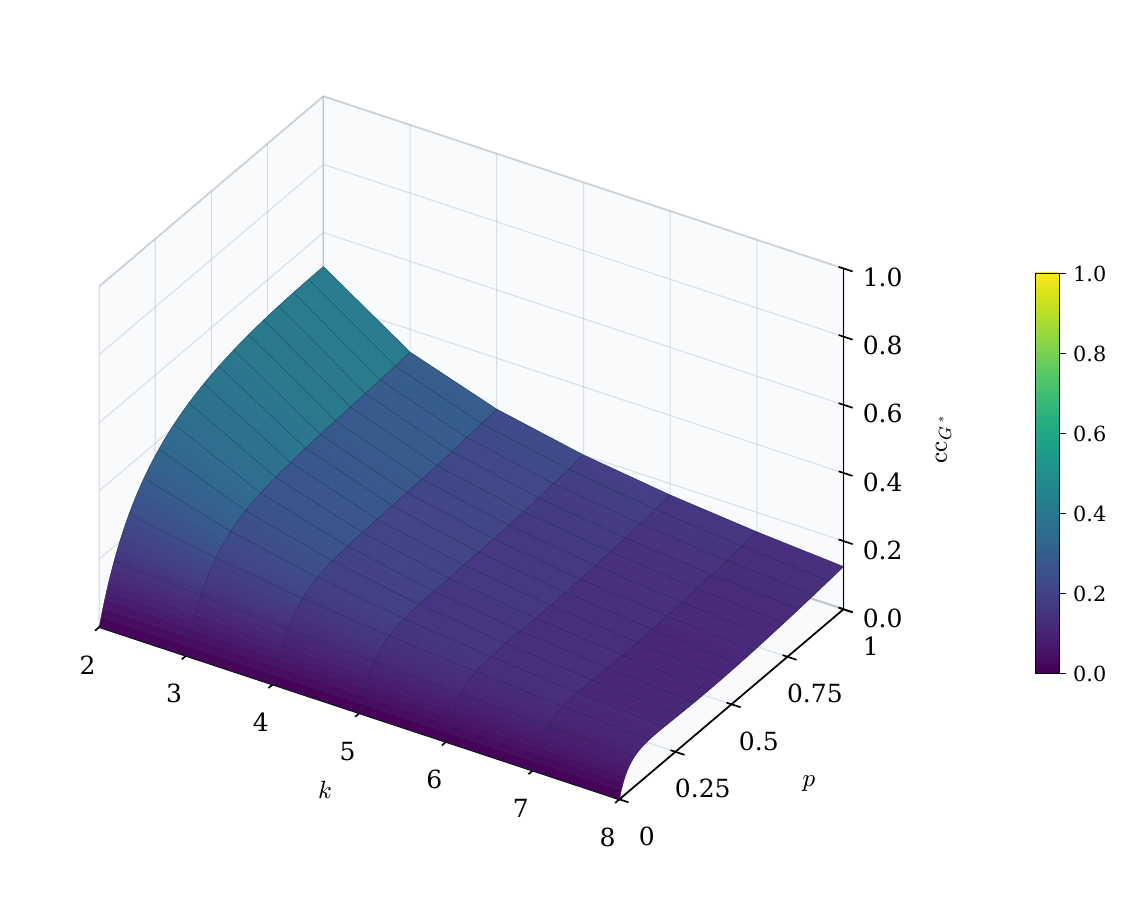}
    \caption{Limiting global clustering coefficient of
    $\cT_p(\mathrm{Regular}(k))$.}
    \label{fig:cc_global_regular_surface}
\end{figure}

In order to apply Corollary~\ref{cor:cc_global}, note that
\[
    \sup_n \bbE \lt d_{H_n}(o_{H_n})^5 \rt = k^5 < \infty.
\]
Thus,
\begin{equation}\label{eq:cc_global_k_regular}
    \begin{split}
    \cc_{G^*} &= \frac{(k (k-1) + 2 k (k-1)) \cdot p + k(k-1)(k-2) \cdot p^3}{k(k-1) + 2 k (k-1) \cdot k \cdot p + k(k-1) \cdot (k(k-1)-1) \cdot p^2}\\
    &= \frac{3p + (k-2)p^3}{1 + 2kp + (k^2 - k -1)p^2}.
    \end{split}
\end{equation}
In Figure~\ref{fig:cc_global_regular_surface}, we plot $\cc_{G^*}$ as a function of $k$ and $p$.

\subsection{Configuration models.} 
Configuration is a generalisation of the random $k$-regular graph. Here for each $n \in \bbN$ let the vertex set of $H_n$ be given by $\cV(H_n) = \{v_i : i \in [n]\}$. We are also given a triangular array of degrees $(m_{n,i})_{i \in [n], n \in \bbN}$ satisfying the conditions that $m_{n,i} \ge 1$ and $\sum_{i=1}^n m_{n,i}$ is even. We shall construct $H_n$ in such a way that $d_{H_n}(v_i) = m_{n,i}$. Note that if sum of the degree sequence is even, then it is always possible to construct a multigraph with the given degree sequence. So, we first construct a random multigraph $\tilde{H}_n$ from the given degree sequence in the following manner. The reader is referred to Chapter~7 of \cite{vanderhofstad2017random} and Chapter~4 of \cite{vanderhofstad2024random} for a comprehensive study of the configuration models.

\subsubsection{The model}
Fix $n \in \bbN$. First endow each vertex $v_i$ with $m_{n,i}$ many half edges. So, there are in total $T_n:= \sum_{i=1}^n m_{n,i}$ many half edges. Take an enumeration of them say $(e_j)_{j \in [T_n]}$. We shall pair up these half edges in an uniform manner to construct an edge in $\tilde{H}_n$. In the first step choose an half edge from $\{e_j: j \ne 1\}$ uniformly and pair it with $e_1$. Inductively, suppose $t$ many edges are already constructed, i.e., $2t$ many half edges are already paired up. In the next step, take the least index $j$ such that $e_j$ is not yet paired up and pair it up uniformly with one of the $T_n-2t-1$ many unpaired half edges. Continue this process until all the half edges are exhausted. Clearly, $\tilde{H}_n$ can be a multigraph, i.e., it can contain self loops and multiedges. Moreover, the degree of the vertex $v_i$ is $m_{n,i}$ where each self loop at $v_i$ contributes $2$ to its degree. The reader can also notice that due to the exchangeable nature of the construction, the distribution of $\tilde{H}_n$ is not affected by the selection of the enumeration of the half edges. 

Once $\tilde{H}_n$ is constructed, the distribution of $H_n$ is given by the distribution of $\tilde{H}_n$ conditional on $\tilde{H}_n$ being simple. In other words, given a simple graph $K$ on the vertex set $\{v_i: i \in [n]\}$,
\[
    \bbP\lt H_n = K \rt = \bbP \lt \tilde{H}_n = K \mid \tilde{H}_n\,\,\text{is simple} \rt.
\]
A caveat is that even when $T_n$ is even, there might not exist any simple graph with the degree sequence $(m_{n,i})_{i \in [n]}$. To resolve this issue and later to apply Theorem~\ref{thm:lwl} we make the following regularity assumptions on the degree sequence.

\begin{assumption}\label{ass:regularity_CM}
    Let $U_n \sim \Unif([n])$. Then there exists a probability mass function $q := (q_i)_{i=1}^\infty$, such that the following holds:
    \begin{align}
        & m_{n, U_n} \convd q, \label{eq:ass_CM_1}\\
        & \bbE \lt m_{n, U_n} \rt \to \sum_{k=1}^\infty k q_k < \infty, \label{eq:ass_CM_2}\\
        & \bbE \lt m^2_{n, U_n} \rt \to \sum_{k=1}^\infty k^2 q_k < \infty. \label{eq:ass_CM_3}
    \end{align}
\end{assumption}
Under the Assumption~\ref{ass:regularity_CM} one can show that 
\[
    \lim_{n \to \infty} \bbP \lt \tilde{H}_n \,\, \text{is simple}\rt = e^{-\frac{\eta}{2} - \frac{\eta^2}{4}} >0
\]
where 
\[
    \eta := \frac{\sum_{k=1}^\infty k(k-1) q_k}{\sum_{k=1}^\infty k q_k} < \infty.
\]
Thus, under the Assumption~\ref{ass:regularity_CM}, $H_n$ is well defined for sufficiently large $n$. 

\subsubsection{Local limit}
The next theorem from \cite{vanderhofstad2024random} describes the local limit of $H_n$.

\begin{theorem}[Local limit of configuration model]\label{thm:lwl_CM}
    Under the Assumption~\ref{ass:regularity_CM}, $H_n$ converges locally in probability to the unimodular Galton-Watson tree with root offspring distribution $(q_k)_{k \in \bbN}$.
\end{theorem}
The unimodular Galton-Watson tree with root offspring distribution $(q_k)_{k \in \bbN}$ is an possibly infinite random tree where the number of offsprings of the root is distributed according to the pmf $(q_k)_{k \in \bbN}$ and the number of offsprings of the vertices in the subsequent generations are distributed according to the size biased distribution of $q$ given by $\lf q^*_k \rf_{k \in \bbN \cup \{0\}}$ where
\[
    q^*_k := \frac{(k+1) q_{k+1}}{\sum_{j=1}^\infty j q_j}\qquad \text{for}\,\, k \in \bbN \cup \{0\}.
\]
Moreover, the number of offsprings for distinct vertices are independent. Now, an application of Theorem~\ref{thm:lwl} yields the local limit of $G_n$.

\subsubsection{Clustering coefficients}

Let us now compute the clustering coefficients. Note that $N \sim (q_k)_{k \in \bbN}$ and conditional on $N$, the variables $d_{H^*}(v_i)-1$, $i \in [N]$, are $\iid$ with pmf $(q^*_k)_{k \in \bbN \cup \{0\}}$, where $v_1, \ldots, v_N$ are the neighbors of $o_{H^*}$ in $H^*$. So,
\begin{align*}
    M &= \sum_{i=1}^N (d_{H^*}(v_i)-1),\\
    P &= \sum_{i=1}^N (d_{H^*}(v_i)-1)(d_{H^*}(v_i)-2).
\end{align*}
This gives the necessary information to compute $\overline{\cc}_{G^*}$ by Corollary~\ref{cor:cc_local}(c).

Now for the convergence of $\cc_{G_n}$, we make the assumption
\[
    \sup_{n} \frac{1}{n} \sum_{i=1}^n m^{4+ \epsilon}_{n,i} < \infty
\]
for some $\epsilon >0$, and $\sum_{k=2}^\infty q_k>0$. By Corollary~\ref{cor:cc_global} and Remark~\ref{rmk:4_moment_UI},
\[
    \cc_{G_n} \convp \cc_{G^*}.
\]
We compute $\cc_{G^*}$ next. To that end, first we fix some notations for the first two moments of $(q_k)_{k \in \bbN}$ and $(q^*_k)_{k \in \bbN \cup \{0\}}$. Define
\begin{align*}
    &\tau_1 := \sum_{k=1}^\infty k q_k,
    & \tau_2 := \sum_{k=1}^\infty k^2 q_k,\\
    & \tau^*_1 := \sum_{k=1}^\infty k q^*_k = \frac{1}{\tau_1} \sum_{k=0}^\infty k (k+1) q_{k+1},
    & \tau^*_2 := \sum_{k=1}^\infty k^2 q^*_k = \frac{1}{\tau_1} \sum_{k=0}^\infty k^2 (k+1) q_{k+1}.
\end{align*}
Now one can easily calculate
\begin{align*}
    &\bbE_\mu[N] = \tau_1,\\
    &\bbE_\mu[N^2] = \tau_2, \\
    &\bbE_\mu[M] = \tau^*_1 \bbE_\mu [N] = \tau_1 \tau^*_1,\\
    &\bbE_\mu [MN] = \tau^*_1 \bbE_\mu[N^2] = \tau_2 \tau^*_1, \\
    &\bbE_\mu[M^2] = {\tau^*_1}^2 \bbE_\mu [N(N-1)] + \tau^*_2 \bbE_\mu[N] = \tau_2 {\tau^*_1}^2 - \tau_1 {\tau^*_1}^2 + \tau_1 \tau^*_2,\\
    &\bbE_\mu[P] = (\tau^*_2 - \tau^*_1) \bbE_\mu[N] = \tau_1 (\tau^*_2 - \tau^*_1).
\end{align*}
Hence
\[
    \cc_{G^*} = \frac{(\tau_2 - \tau_1 + 2 \tau_1 \tau^*_1) p + \tau_1 (\tau^*_2 - \tau^*_1) p^3}{(\tau_2 - \tau_1) + 2 \tau_2 \tau_1^* p + (\tau_2 {\tau^*_1}^2 - \tau_1 {\tau^*_1}^2 + \tau_1 \tau^*_2 - \tau_1 \tau^*_1) p^2}.
\]

\subsection{Preferential attachment models.} \label{subsec:PA}

The preferential attachment model was first introduced by Barab{\'a}si and Albert in their seminal work \cite{barabasi1999emergence}. We consider a more general version of the model introduced in \cite{berger2014asymptotic}. In \cite{berger2014asymptotic}, the authors work out the local weak limit of this model. We introduce the model now.

\subsubsection{The model}
Fix $\alpha \in [0,1)$ and $m \ge 2$. We define $(\tilde H_n)_{n \ge 1}$, a sequence of (unrooted) dynamic undirected multigraphs without self loops, inductively as follows. Let $\tilde H_1$ be the graph with one vertex $v_1$ and no edge and $\tilde H_2$ be the multigraph with two vertices $v_1, v_2$ and $m$ edges between them. At each subsequent stage, $m$ edges are added. So, for $n \ge 1$, $\tilde H_n$ has $(n-1)m$ edges, counted with multiplicity. Now, suppose for some $n \ge 3$, $\tilde H_{n-1}$ is already constructed. Define $\cV(\tilde H_n) := \cV(\tilde H_{n-1}) \cup \{v_n\}$. Conditional on $\tilde H_{n-1}$, we add $m$ $\iid$ edges $e_1, \ldots, e_m$ with distribution
\[
    \bbP \lt e_1 = \{v_n, v_i\} \mid \tilde H_{n-1} \rt
    = \alpha \frac{1}{n-1} + (1 - \alpha) \frac{d_{\tilde H_{n-1}}(v_i)}{2(n-2)m}
    \qquad \text{for}\,\, i \in [n-1],
\]
where degrees are counted with multiplicity. Then, $(\tilde H_n)_{n \in \bbN}$ is called an \emph{independent preferential attachment model} with parameters $\alpha$ and $m$.

Let $H_n$ be the underlying simple graph of $\tilde H_n$, obtained by replacing each set of multiedges by a single edge. The attachment probabilities above always use the degrees in $\tilde H_n$. We apply triadic closure to $H_n$ and write $G_n := \cT_p(H_n)$, where $p \in [0,1]$. The local convergence results in \cite{berger2014asymptotic,vanderhofstad2024random} give the \emph{P\'{o}lya point tree}, denoted by $H^* = \PPT(\alpha,m)$, as the local limit in probability. Note that erasing parallel edges preserves this limit. Indeed by Theorem~2.15 of \cite{vanderhofstad2024random}, it is enough to verify \eqref{eq:lwc_in_prob_equivalent} only by choosing $K^*$ to be simple rooted trees, and in that case for any $v \in \cV(H_n)$, $r \ge 0$, if $B_{\tilde{H}_n}(v, r) \cong K^*$, then $B_{H_n}(v, r) \cong K^*$, implying, in particular, that
\[
    p_{\tilde{H}_n}(K^*, r) \le p_{H_n}(K^*, r).
\]
Moreover,
\begin{align*}
    0 \le p_{H_n}(K^*, r) - p_{\tilde{H}_n}(K^*, r) &\le \frac{1}{n} \sum_{v \in \cV(\tilde{H}_n)} \bbE [\bbI( B_{\tilde{H}_n} (v, r) \,\,\text{is not simple})] \\
    & \convp \mu [B_{H^*} (o_{H^*}, r) \,\,\text{is not simple}] = 0,
\end{align*}
Where the convergence in the second line follows from the definition of local convergence \ref{def:lwc}. Thus, $H_n$ converges locally in probability to $H^*$. Write $\mu$ for the law of $H^*$ and $G^* := \cT_p(\PPT(\alpha,m))$. Computation of $\overline{\cc}_{G^*}$ and $\cc_{G^*}$ requires the knowledge of the $2$-neighborhood of $\PPT(\alpha,m)$. We describe it in detail now. For a complete description of $\PPT(\alpha,m)$ the reader is directed to \cite{berger2014asymptotic,vanderhofstad2024random}.

\subsubsection{Description of $\PPT(\alpha, m)$}
The root is denoted by $(0)$. The neighbors $(0,i)$, $1 \le i \le m$, are older than the root, and the remaining neighbors are younger than the root. The vertex set of the $1$-neighborhood of $(0)$ is given by
\[
    \cV(B_{H^*}((0), 1)) := \{(0)\} \cup \{ (0,i): i \in [m+Q_{(0)}] \}.
\]
Now, for $1 \le i \le m$,
\[
    \cV(B_{H^*}((0,i), 1)) := \{(0), (0,i)\} \cup \{ (0, i, j): j \in [m+Q_{(0,i)}] \},
\]    
and for $m+1 \le i \le m+ Q_{(0)}$,
\[
    \cV(B_{H^*}((0,i), 1)) := \{(0), (0,i)\} \cup \{ (0, i, j): j \in [m -1 + Q_{(0,i)}] \}.
\]
Note that 
\begin{align*}
    d_{H^*}((0)) &= m + Q_{(0)},\\
    d_{H^*}((0,i)) &= \begin{cases}
        m +1 + Q_{(0,i)} \qquad \text{if}\,\, 1 \le i \le m, \\
        m+ Q_{(0,i)} \qquad \text{if}\,\, m+1 \le i \le m+ Q_{(0)}.
    \end{cases}
\end{align*}
Then,
\begin{align*}
    & N = m + Q_{(0)}, \\
    & M = m^2 + (m-1) Q_{(0)} + \sum_{i=1}^m Q_{(0,i)} + \sum_{i=m+1}^{m+Q_{(0)}} Q_{(0,i)}, \\
    & P = \sum_{i=1}^m (m+Q_{(0,i)})(m-1+Q_{(0,i)}) \\
    &\qquad\quad + \sum_{i=m+1}^{m+Q_{(0)}} (m-1+Q_{(0,i)})(m-2+Q_{(0,i)}).
\end{align*}
We now describe the distributions of $Q_{(0)}$ and $Q_{(0,i)}$'s. There is a collection of positions $(x_v)_{v \in \cV(H^*)}$ in $[0,1]$ associated with the \emph{P\'{o}lya point process}. Together with the Gamma variables below, these positions determine the conditional distributions of the numbers of younger neighbors. We only describe the distribution of $Q_{(0)}$ and $\lf Q_{(0,i)} \rf_{i \in [N]}$ here. Suppose  
\begin{align*}
    U_{(0)}, U_{(0,1)}, \ldots, U_{(0,m)} &\overset{\iid}{\sim} \Unif([0,1]), \\
    T_{(0)}, T_{(0, m+1)}, T_{(0, m+2)}, \ldots &\overset{\iid}{\sim} \mathrm{Gamma} \left(\frac{(1 + \alpha)m}{1- \alpha}, 1 \right), \\
    T_{(0,1)}, \ldots, T_{(0,m)} &\overset{\iid}{\sim} \mathrm{Gamma}\left( \frac{(1 + \alpha)m}{1- \alpha} + 1, 1 \right),
\end{align*}
are independent collections, where the Gamma distributions have rate $1$. Firstly, define 
\begin{align*}
    x_{(0)} &:= U_{(0)}^{\frac{1 + \alpha}{2}},\\
    x_{(0, i)} &:= x_{(0)} U_{(0,i)} \qquad \text{for}\,\, 1 \le i \le m.
\end{align*}
Now, conditional on $x_{(0)}$ and $T_{(0)}$,
\[
    \{ x_{(0,m+1)}, \ldots, x_{(0, m + Q_{(0)})} \} \sim \PPP(\rho_{(0)}),
\]
where $\PPP$ stands for \emph{Poisson point process} and 
\[
    \rho_{(0)}(x) := \frac{1 - \alpha}{1 + \alpha} T_{(0)} x_{(0)}^{-\frac{1 - \alpha}{1 + \alpha}} x^{- \frac{2 \alpha}{1+\alpha}} \bbI_{[x_{(0)}, 1]}(x). 
\]
Attach the independent Gamma marks $T_{(0,i)}$ to the younger-neighbor points in increasing order of their positions. Lastly, conditional on the construction so far, define independent point processes: 
\begin{align*}
    & \{ x_{(0,i,m+1)}, \ldots, x_{(0,i, m+ Q_{(0,i)})} \} \sim \PPP(\rho_{(0, i)}) \qquad \text{if}\,\, 1 \le i \le m,\\
    & \{ x_{(0,i,m)}, \ldots, x_{(0,i, m -1 + Q_{(0,i)})} \} \sim \PPP(\rho_{(0, i)}) \qquad \text{if}\,\, m + 1 \le i \le m + Q_{(0)},
\end{align*}
where
\[
    \rho_{(0,i)}(x) := \frac{1 - \alpha}{1 + \alpha} T_{(0, i)} x_{(0,i)}^{-\frac{1 - \alpha}{1 + \alpha}} x^{- \frac{2 \alpha}{1+\alpha}} \bbI_{[x_{(0,i)}, 1]}(x).
\]
Thus, conditional on $x_{(0)}$, $T_{(0)}$,
\begin{equation}\label{eq:P_0}
    Q_{(0)} \sim \Poisson \left(\lambda_{(0)} \right),
\end{equation}
where
\begin{align*}
    \lambda_{(0)} &:= \int_{x_{(0)}}^1 \rho_{(0)} (x) \, dx \\
    &= \frac{1 - \alpha}{1 + \alpha} T_{(0)} x_{(0)}^{-\frac{1 - \alpha}{1 + \alpha}} \int_{x_{(0)}}^1 x^{- \frac{2 \alpha}{1+\alpha}} \, dx \\
    &= T_{(0)} \left( x_{(0)}^{-\frac{1 - \alpha}{1 + \alpha}} - 1 \right) \\
    &= T_{(0)} \left( U_{(0)}^{-\frac{1 - \alpha}{2}} - 1 \right).
\end{align*}
Next conditional on $\left( x_{(0,i)} \right)_{i \in [m + Q_{(0)}]}$ and $\left( T_{(0,i)} \right)_{i \in [m + Q_{(0)}]}$, for $i \in [m + Q_{(0)}]$ we have the independent variables
\begin{align}
    Q_{(0,i)}  \sim \Poisson(\lambda_{(0,i)}), \label{eq:P_(0,i)} 
\end{align}
where
\begin{align*}
    \lambda_{(0,i)} &:= \int_{x_{(0,i)}}^1 \rho_{(0,i)} (x) \, dx \\
    &= \frac{1 - \alpha}{1 + \alpha} T_{(0, i)} x_{(0,i)}^{-\frac{1 - \alpha}{1 + \alpha}} \int_{x_{(0,i)}}^1 x^{- \frac{2 \alpha}{1+\alpha}} \, dx \\
    &= T_{(0,i)} \left( x_{(0,i)}^{-\frac{1- \alpha}{1 + \alpha}} -1 \right).
\end{align*}
For an older neighbor, $Q_{(0,i)}$ counts its additional younger neighbors, excluding the root. On the other hand, for a younger neighbor, it counts all its younger neighbors. This gives the joint distribution of $M$, $N$ and $P$. In particular, by Corollary~\ref{cor:cc_local}(c), for every $\alpha \in [0,1)$,
\[
    \overline{\cc}_{G_n} \convp \overline{\cc}_{G^*}
    = \bbE_\mu \left[
        \frac{N(N-1)p+2X+P\cdot\frac{X(X-1)}{M(M-1)}p}
        {N(N-1)+(2N-1)X+X^2}
      \right],
\]
where, conditional on $H^*$, $X \sim \Binom(M,p)$, and $0/0$ is interpreted as $0$. The expectation includes this additional binomial randomness. Next we verify uniform integrability and compute $\cc_{G^*}$. 

\subsubsection{Uniform integrability of $\left(|\cV(B_{H_n}(o_{H_n},2))|^2\right)_{n\in\bbN}$}

We first show that if $\alpha \in (\frac12,1)$, then
\[
    \sup_n \bbE \lt d_{H_n}(o_{H_n})^{4+\epsilon} \rt < \infty
\]
for some $\epsilon>0$. In general, we show the following.
\begin{lemma}\label{lem:degree_moment_bound_PA}
    For every real $r\ge1$ and $1\le i\le n$, we have
    \begin{equation}\label{eq:degree_PA_bound}
        \bbE \lt d_{H_n}(v_i)^r \rt
        \le \bbE \lt d_{\tilde H_n}(v_i)^r \rt
        \le C_{m,\alpha,r} n^{\frac{r(1-\alpha)}2} i^{-\frac{r(1-\alpha)}2},
    \end{equation}
    for some positive constant $C_{m,\alpha,r}$ depending only on $m$, $\alpha$ and $r$.
\end{lemma}

\begin{proof}[Proof of Lemma~\ref{lem:degree_moment_bound_PA}]
    Since $d_{H_n}(v_i)\le d_{\tilde H_n}(v_i)$, it is enough to prove the second inequality. Put $\beta=(1-\alpha)/2>0$. We first prove the bound for integer $r$ by induction. The case $n=1$ is immediate. For $n\ge2$, fix $i\le n$ and put $j=\max\{i,2\}$. Then $d_{\tilde H_j}(v_i)=m$. For $k\ge j$,
    \[
        d_{\tilde H_{k+1}}(v_i)=d_{\tilde H_k}(v_i)+\chi,
    \]
    where, conditional on $\tilde H_k$,
    \[
        \chi\sim\Binom\lf m,\frac{\alpha}{k}
              +\frac{(1-\alpha)d_{\tilde H_k}(v_i)}{2m(k-1)}\rf.
    \]
    Here $\chi$ denotes the degree increment at step $k$; its distribution depends on $k$ and $i$. Thus,
    \begin{equation}\label{eq:chi_expec}
        \bbE\lt\chi\mid\tilde H_k\rt
        =\frac{m\alpha}{k}+\frac{\beta d_{\tilde H_k}(v_i)}{k-1}.
    \end{equation}
    We shall use, for fixed $c>0$ and $2\le a\le b$,
    \[
        \prod_{k=a}^b\lf1+\frac{c}{k-1}\rf
        \le C_c\lf\frac{b+1}{a}\rf^c,
    \]
    which follows by taking logarithms and using $\log(1+x)\le x$. Empty products equal $1$, and the same bound holds when $a=b+1$, on taking $C_c\ge1$.

    For $r=1$, iteration gives
    \begin{align*}
        \bbE\lt d_{\tilde H_n}(v_i)\rt
        &=m\prod_{k=j}^{n-1}\lf1+\frac{\beta}{k-1}\rf
          +m\alpha\sum_{k=j}^{n-1}\frac1k
            \prod_{l=k+1}^{n-1}\lf1+\frac{\beta}{l-1}\rf\\
        &\le C_{m,\alpha,1}n^\beta
             \lf j^{-\beta}+\sum_{k=j}^{n-1}k^{-1-\beta}\rf\\
        &\le C'_{m,\alpha,1}n^\beta i^{-\beta}.
    \end{align*}
    Here we used
    \[
        \sum_{k=j}^{n-1}k^{-1-\beta}
        \le j^{-1-\beta}+\int_j^\infty x^{-1-\beta}\,dx
        \le \lf1+\frac1\beta\rf j^{-\beta},
    \]
    and $j\ge i$. This proves \eqref{eq:degree_PA_bound} for $r=1$.

    Now suppose the bound holds for $r-1$, where $r\ge2$ is an integer. To estimate the increment, we use the following consequence of Taylor's theorem: for every real $r\ge1$ and $x\in[0,1]$,
    \[
        (1+x)^r\le1+rx+C_r x^2,
        \qquad
        C_r:=\frac{r(r-1)}2\max\{1,2^{r-2}\}.
    \]
    Indeed, the second derivative of $(1+x)^r$ on $[0,1]$ is bounded above by $2C_r$. Since $0\le\chi\le m\le d_{\tilde H_k}(v_i)$, we have
    \[
        0\le\frac{\chi}{d_{\tilde H_k}(v_i)}\le1,
        \qquad
        \chi^2\le m\chi.
    \]
    Applying the preceding inequality and then \eqref{eq:chi_expec}, we obtain
    \begin{equation}\label{eq:binomial_thm_bound}
        \begin{split}
            &\bbE\lt d_{\tilde H_{k+1}}(v_i)^r\mid\tilde H_k\rt\\
            &\le d_{\tilde H_k}(v_i)^r
            +r d_{\tilde H_k}(v_i)^{r-1} \bbE\lt\chi\mid\tilde H_k\rt\\
            &\qquad
            +C_r d_{\tilde H_k}(v_i)^{r-2}
            \bbE\lt\chi^2\mid\tilde H_k\rt\\
            &\le d_{\tilde H_k}(v_i)^r
            +\lf r d_{\tilde H_k}(v_i)^{r-1}
            +mC_r d_{\tilde H_k}(v_i)^{r-2}\rf
            \lf\frac{m\alpha}{k}
            +\frac{\beta d_{\tilde H_k}(v_i)}{k-1}\rf\\
            &\le\lf1+\frac{r\beta}{k-1}\rf d_{\tilde H_k}(v_i)^r
            +\frac{m\{r\alpha+C_r(\alpha+\beta)\}}{k-1}
            d_{\tilde H_k}(v_i)^{r-1}\\
            &\le\lf1+\frac{r\beta}{k-1}\rf d_{\tilde H_k}(v_i)^r
            + \frac{C'_{m, \alpha, r}}{k-1}
            d_{\tilde H_k}(v_i)^{r-1}.
        \end{split}
    \end{equation}
    In the last step, we used $k^{-1}\le(k-1)^{-1}$ and
    $m d_{\tilde H_k}(v_i)^{r-2}\le d_{\tilde H_k}(v_i)^{r-1}$.
    In particular, all terms other than the linear preferential-attachment term are bounded by a constant times $(k-1)^{-1}d_{\tilde H_k}(v_i)^{r-1}$.

    Taking expectations and iterating gives
    \begin{align*}
        \bbE\lt d_{\tilde H_n}(v_i)^r\rt
        &\le m^r\prod_{k=j}^{n-1}\lf1+\frac{r\beta}{k-1}\rf\\
        &\quad+C'_{m,\alpha,r}\sum_{k=j}^{n-1}
        \frac{\bbE\lt d_{\tilde H_k}(v_i)^{r-1}\rt}{k-1} \prod_{l=k+1}^{n-1}\lf1+\frac{r\beta}{l-1}\rf\\
        &\le C''_{m,\alpha,r}n^{r\beta}j^{-r\beta}\\
        &\quad+C'''_{m,\alpha,r}n^{r\beta}
        \sum_{k=j}^{n-1}k^{-1-r\beta}
        \bbE\lt d_{\tilde H_k}(v_i)^{r-1}\rt\\
        &\le C''_{m,\alpha,r}n^{r\beta}i^{-r\beta}
        +C''''_{m,\alpha,r}n^{r\beta}i^{-(r-1)\beta}
        \sum_{k=j}^{n-1}k^{-1-\beta}\\
        &\le C_{m,\alpha,r}n^{r\beta}i^{-r\beta},
    \end{align*}
    where the penultimate inequality follows from the induction hypothesis. All constants here depend only on $m$, $\alpha$ and $r$. This completes the induction.

    For real $r\ge1$, put $s=\lceil r\rceil$. By the integer-moment bound and Jensen's inequality,
    \[
        \bbE\lt d_{\tilde H_n}(v_i)^r\rt
        \le\lf\bbE\lt d_{\tilde H_n}(v_i)^s\rt\rf^{r/s}
        \le C_{m,\alpha,s}^{r/s}(n/i)^{r\beta}.
    \]
    This proves \eqref{eq:degree_PA_bound} for every real $r\ge1$.
\end{proof}

In particular, averaging \eqref{eq:degree_PA_bound} over $i$ yields
\begin{align*}
    \bbE\lt d_{H_n}(o_{H_n})^r\rt
    &=\frac1n\sum_{i=1}^n\bbE\lt d_{H_n}(v_i)^r\rt\\
    &\le C_{m,\alpha,r}n^{r\beta-1}\sum_{i=1}^n i^{-r\beta}\\
    &\le \frac{C_{m,\alpha,r}}{1-r\beta},
\end{align*}
whenever $r\beta=r(1-\alpha)/2<1$. If $\alpha>1/2$, choose
\[
    0<\epsilon<\frac{2}{1-\alpha}-4.
\]
Then $(4+\epsilon)(1-\alpha)/2<1$, and hence
\[
    \sup_n\bbE\lt d_{H_n}(o_{H_n})^{4+\epsilon}\rt<\infty.
\]
So, $\lf d_{H_n}(o_{H_n})^{4} \rf_{n \in \bbN}$ is uniformly integrable. By Remark~\ref{rmk:4_moment_UI}, $\left(|\cV(B_{H_n}(o_{H_n},2))|^2\right)_{n\in\bbN}$ is also uniformly integrable.

\begin{remark}
    The degree-moment condition is sharp. Indeed, by Lemma~5.2 of \cite{berger2014asymptotic},
    \[
        \bbP_\mu\lt d_{H^*}(o_{H^*})=m+k\rt
        =C k^{-\frac{1+\alpha}{1-\alpha}-2}
             \lf1+O\lf\frac1k\rf\rf.
    \]
    Hence, for every real $r\ge1$, $\bbE_\mu\lt d_{H^*}(o_{H^*})^r\rt<\infty$ if and only if $r(1-\alpha)/2<1$. For each fixed $K>0$, local convergence gives
    \[
        \bbE_\mu\lt\min\{d_{H^*}(o_{H^*})^r,K\}\rt
        =\lim_{n\to\infty}
          \bbE\lt\min\{d_{H_n}(o_{H_n})^r,K\}\rt
        \le\sup_n\bbE\lt d_{H_n}(o_{H_n})^r\rt.
    \]
    Letting $K\to\infty$ and using the preceding moment bound, we therefore obtain
    \[
        \sup_n\bbE\lt d_{H_n}(o_{H_n})^r\rt<\infty
        \quad\Longleftrightarrow\quad \alpha>1-\frac2r.
    \]
    Thus $\alpha>1/2$ is the best range obtainable from the degree-moment criterion in Remark~\ref{rmk:4_moment_UI}. 
\end{remark}

\subsubsection{Calculation of $\cc_{G^*}$}

In the rest of this section we assume that $\alpha\in(\frac12,1)$, so that the uniform integrability established above applies. The expectation identities below in fact hold for $\alpha>1/3$. We write $\bbE_\mu$ also for expectations over the auxiliary variables in the construction of $\PPT(\alpha,m)$. The following computations require various expectation calculations. These are collected in Lemma~\ref{lem:expectation_computations_PA}, whose proof we defer to Appendix~\ref{sec:aux}.

\begin{figure}[htbp]
    \centering
    \begin{minipage}[t]{0.49\linewidth}
        \centering
        \includegraphics[width= 0.8\linewidth]{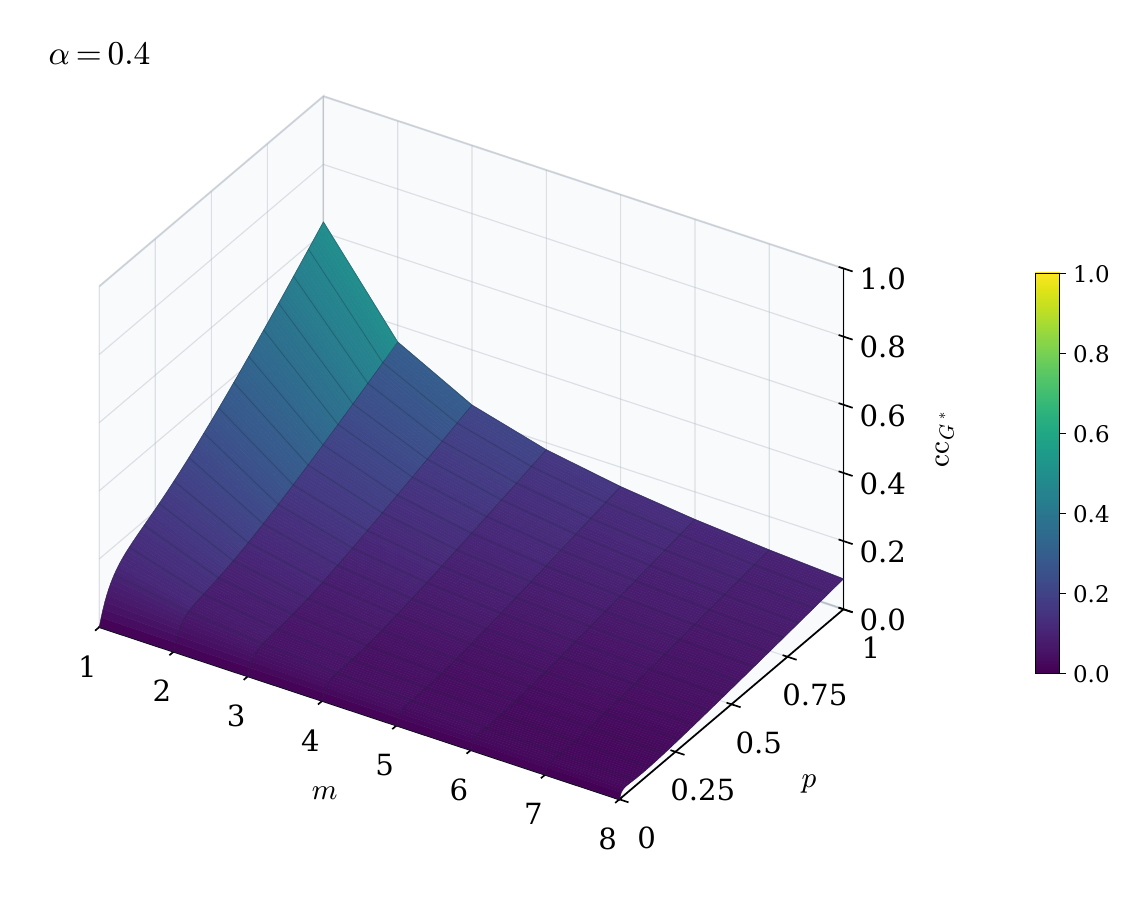}
        \par\smallskip
        {\small (a) $\alpha=0.4$}
    \end{minipage}\hfill
    \begin{minipage}[t]{0.49\linewidth}
        \centering
        \includegraphics[width= 0.8\linewidth]{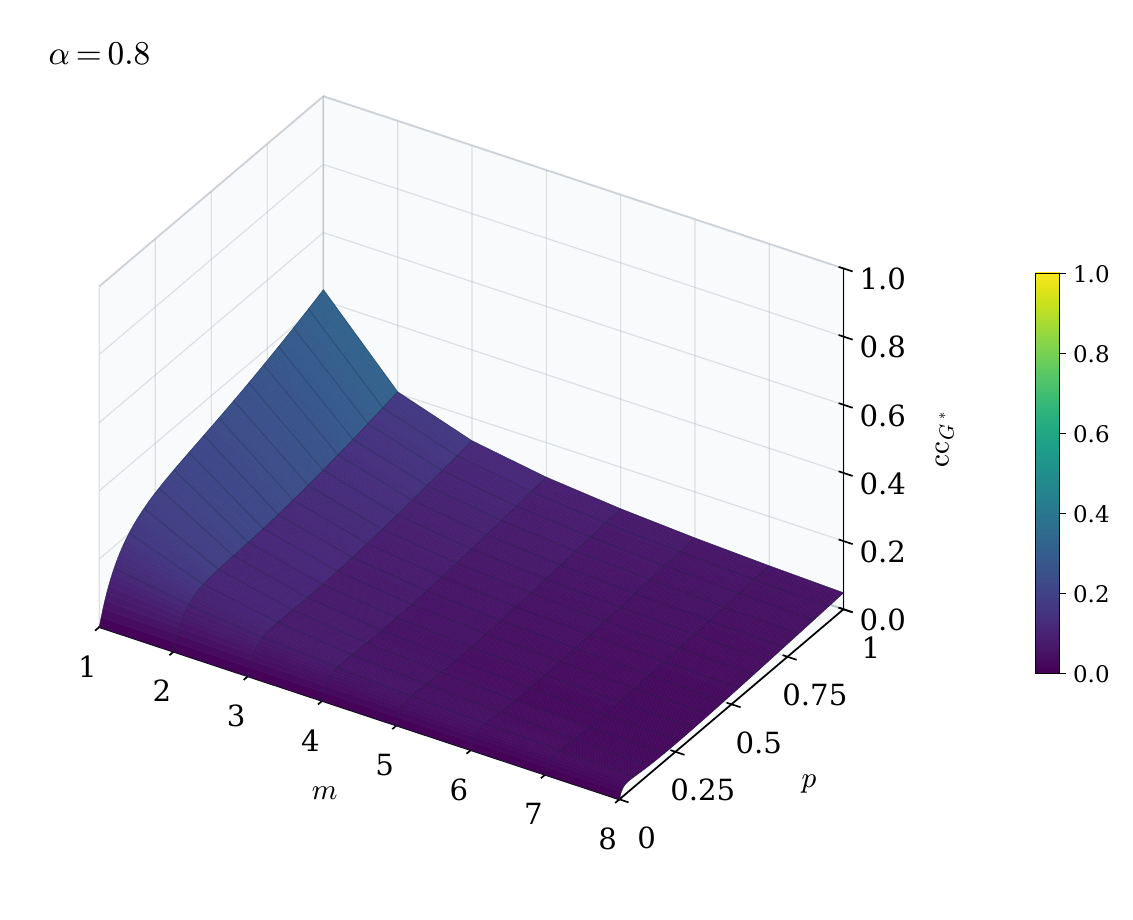}
        \par\smallskip
        {\small (b) $\alpha=0.8$}
    \end{minipage}
    \caption{Limiting global clustering coefficient $\mathrm{cc}_{G^*}$ of $\cT_p(\mathrm{PA}(\alpha, m))$}
    \label{fig:pa_global_clustering_surfaces}
\end{figure}

\begin{lemma}\label{lem:expectation_computations_PA}
For $\alpha>1/3$, we have the following identities.
\begin{align}
    & \bbE_\mu \lt \lambda_{(0)} \rt = m, \label{eq:estimate_1}\\
    & \bbE_\mu \lt \lambda_{(0)}^2 \rt = \frac{m ((1 + \alpha) m + (1 - \alpha))}{ \alpha} , \label{eq:estimate_2}\\
    & \bbE_\mu \lt \lambda_{(0,1)} \rt = \frac{(1 + \alpha) m + (1 - \alpha)}{ \alpha}, \label{eq:estimate_3} \\
    & \bbE_\mu \lt \lambda_{(0)} \lambda_{(0,1)}\rt = \frac{m (1 + 2 \alpha) ((1 + \alpha)m + (1 - \alpha))}{2 \alpha^2}, \label{eq:estimate_4} \\
    &\bbE_\mu\lt\lambda_{(0,1)}^2\rt
    =\frac{3((1+\alpha)m+(1-\alpha))((1+\alpha)m+2(1-\alpha))}
              {\alpha(3\alpha-1)}, \label{eq:estimate_4.5}\\
    & \bbE_\mu \lt \lf  \sum_{i=1}^m Q_{(0,i)}\rf^2 \rt = \frac{
    \begin{aligned}
        &m ((1 + \alpha)m + (1 - \alpha)) ((-1 - 2 \alpha + 11 \alpha^2 + 12 \alpha^3) m^2 \\
        &\qquad \qquad \qquad \qquad +(2 \alpha + 14 \alpha^2 - 12 \alpha^3) m + (1 + 7 \alpha^2))
    \end{aligned}
    }{4 \alpha^3 (-1 + 3 \alpha)}, \label{eq:estimate_5} \\
    & \bbE_\mu \lt \sum_{i=m+1}^{m+ Q_{(0)}} \lambda_{(0,i)} \rt = m^2, \label{eq:estimate_6} \\
    & \bbE_\mu \lt \lambda_{(0,1)} \sum_{i= m +1}^{m + Q_{(0)}} \lambda_{(0,i)} \rt = \frac{m^2 ((1+\alpha)m+(1-\alpha))(1+3\alpha+4\alpha^2)}{4\alpha^3}, \label{eq:estimate_7} \\
    &\bbE_\mu \lt Q_{(0)} \sum_{i= m +1}^{m + Q_{(0)}} \lambda_{(0,i)} \rt = \frac{m^2 (1 + \alpha) ( (1 + 2 \alpha) m +1)}{2 \alpha^2} , \label{eq:estimate_8} \\
    &\bbE_\mu\lt\sum_{i=m+1}^{m+Q_{(0)}}\lambda_{(0,i)}^2\rt
    =\frac{m^2((1+\alpha)m+(1-\alpha))}{\alpha}, \label{eq:estimate_8.5}\\
    & \bbE_\mu \lt \lf \sum_{i=m+1}^{m+Q_{(0)}} Q_{(0,i)}\rf^2 \rt
    =m^2+\frac{m^2((1+\alpha)m+(1-\alpha))}{\alpha}\nonumber\\
    &\qquad\qquad\qquad\qquad
      +\frac{m^3(1+\alpha)(1+2\alpha)((1+\alpha)m+(1-\alpha))}{2\alpha^3}.
      \label{eq:estimate_9}
\end{align}
\end{lemma}
Now, using Lemma~\ref{lem:expectation_computations_PA}, we have the following identities.
\begin{align*}
    \bbE_\mu [N] &= m + \bbE_\mu \lt Q_{(0)}\rt \\
    &= m + \bbE_\mu \lt \lambda_{(0)} \rt \\
    &= m + m \\
    &= 2m,\\
    \bbE_\mu [N^2] & = m^2 + 2m \bbE_\mu \lt Q_{(0)} \rt + \bbE_\mu \lt Q_{(0)}^2 \rt \\
    &= m^2 + (2m +1) \bbE_\mu \lt \lambda_{(0)} \rt + \bbE_\mu \lt \lambda_{(0)}^2 \rt \\
    &= m^2 + (2m +1) m + \frac{m ((1 + \alpha) m + (1 - \alpha))}{\alpha} \\
    &= \lf 4 + \frac{1}{\alpha} \rf m^2 + \frac{1}{\alpha} m , \\
    \bbE_\mu [M] &= m^2 + (m-1) \bbE_\mu\lt Q_{(0)}\rt + \bbE_\mu \lt \sum_{i=1}^{m} Q_{(0,i)} \rt + \bbE_\mu \lt \sum_{i=m+1}^{m+ Q_{(0)}} Q_{(0,i)} \rt \\
    &= m^2 + (m-1) \bbE_\mu \lt \lambda_{(0)} \rt + m \bbE_\mu \lt \lambda_{(0,1)} \rt + \bbE_\mu \lt \sum_{i=m+1}^{m+ Q_{(0)}} \lambda_{(0,i)} \rt \\
    &= \lf 4 + \frac{1}{\alpha} \rf m^2 + \lf \frac{1}{\alpha} - 2 \rf m ,\\
    \bbE_\mu [MN] &= m^3 + (2 m^2 - m) \bbE_\mu \lt Q_{(0)} \rt + (m-1) \bbE_\mu \lt Q_{(0)}^2 \rt + m \bbE_\mu \lt \sum_{i=1}^m Q_{(0,i)} \rt + \bbE_\mu \lt Q_{(0)} \sum_{i=1}^m Q_{(0,i)} \rt \\
    &+ m \bbE_\mu \lt \sum_{i=m+1}^{m+ Q_{(0)}} Q_{(0,i)} \rt + \bbE_\mu \lt Q_{(0)} \sum_{i=m+1}^{m+ Q_{(0)}} Q_{(0,i)} \rt \\
    &= m^3 + (2 m^2 - 1) \bbE_\mu \lt \lambda_{(0)} \rt + (m-1) \bbE_\mu \lt \lambda_{(0)}^2 \rt + m^2 \bbE_\mu \lt \lambda_{(0,1)} \rt + m \bbE_\mu \lt \lambda_{(0)} \lambda_{(0,1)} \rt \\
    &+ m \bbE_\mu \lt \sum_{i=m+1}^{m+ Q_{(0)}} \lambda_{(0,i)} \rt + \bbE_\mu \lt Q_{(0)} \sum_{i=m+1}^{m+ Q_{(0)}} \lambda_{(0,i)} \rt \\
    &= m^3 + (2 m^2 - 1) m + (m-1) \frac{m ((1 + \alpha) m + (1 - \alpha))}{\alpha} + m^2 \frac{(1 + \alpha) m + (1 - \alpha)}{ \alpha} \\
    & + m \frac{m (1 + 2 \alpha) ((1 + \alpha)m + (1 - \alpha))}{2 \alpha^2} + m m^2 + \frac{m^2 (1 + \alpha) ( (1 + 2 \alpha) m +1)}{2 \alpha^2} \\
    &= \lf 8 + \frac{5}{\alpha} + \frac{1}{\alpha^2} \rf m^3 + \lf -4 + \frac{2}{\alpha} + \frac{1}{\alpha^2} \rf m^2 - \frac{1}{\alpha} m, 
\end{align*}
\begin{align*}
    &\bbE_\mu [M^2] \\
    &= m^4 + (m-1)^2\bbE_\mu \lt Q_{(0)}^2 \rt + \bbE_\mu \lt \lf  \sum_{i=1}^m Q_{(0,i)}\rf^2 \rt + \bbE_\mu \lt \lf  \sum_{i= m +1}^{m + Q_{(0)}} Q_{(0,i)}\rf^2 \rt \\
    &+ 2 m^2 (m-1) \bbE_\mu \lt Q_{(0)} \rt + 2 m^2 \bbE_\mu \lt \sum_{i=1}^m Q_{(0,i)} \rt + 2 m^2 \bbE_\mu \lt \sum_{i= m +1}^{m + Q_{(0)}} Q_{(0,i)} \rt \\
    &+ 2(m-1) \bbE_\mu \lt Q_{(0)} \sum_{i=1}^m Q_{(0,i)} \rt + 2 (m-1) \bbE_\mu \lt Q_{(0)} \sum_{i= m +1}^{m + Q_{(0)}} Q_{(0,i)} \rt \\
    &+ 2 \bbE_\mu \lt \lf \sum_{i= 1 }^m Q_{(0,i)} \rf \lf \sum_{i= m +1}^{m + Q_{(0)}} Q_{(0,i)} \rf \rt \\
    &= m^4 + (2 m^3 - m^2 -2m +1) \bbE_\mu \lt \lambda_{(0)} \rt + (m-1)^2\bbE_\mu \lt \lambda_{(0)}^2 \rt  + \bbE_\mu \lt \lf  \sum_{i=1}^m Q_{(0,i)}\rf^2 \rt \\
    &+ \bbE_\mu \lt \lf  \sum_{i= m +1}^{m + Q_{(0)}} Q_{(0,i)}\rf^2 \rt + 2 m^3 \bbE_\mu \lt \lambda_{(0,1)} \rt + 2 m^2 \bbE_\mu \lt \sum_{i= m +1}^{m + Q_{(0)}} \lambda_{(0,i)} \rt \\
    &+ 2m(m-1) \bbE_\mu \lt \lambda_{(0)} \lambda_{(0,1)} \rt + 2 (m-1) \bbE_\mu \lt Q_{(0)} \sum_{i= m +1}^{m + Q_{(0)}} \lambda_{(0,i)} \rt + 2 m \bbE_\mu \lt \lambda_{(0,1)} \sum_{i= m +1}^{m + Q_{(0)}} \lambda_{(0,i)} \rt \\
    &= \frac{5+30\alpha+69\alpha^2+64\alpha^3}{4\alpha^3}m^4\\
    &\quad+\frac{-5+5\alpha+61\alpha^2+35\alpha^3-192\alpha^4}
                   {4\alpha^3(3\alpha-1)}m^3\\
    &\quad+\frac{1+11\alpha+3\alpha^2-55\alpha^3+48\alpha^4}
                   {4\alpha^3(3\alpha-1)}m^2
          +\frac{1-\alpha+3\alpha^2+5\alpha^3}
                   {4\alpha^3(3\alpha-1)}m.
\end{align*}

\begin{align*}
    \bbE_\mu[P]
    &=m^2(m-1)+2m^2\bbE_\mu\lt\lambda_{(0,1)}\rt
    +m\bbE_\mu\lt\lambda_{(0,1)}^2\rt\\
    &\quad+(m-1)(m-2)\bbE_\mu\lt Q_{(0)}\rt
    +2(m-1)\bbE_\mu\lt\sum_{i=m+1}^{m+Q_{(0)}}\lambda_{(0,i)}\rt
    +\bbE_\mu\lt\sum_{i=m+1}^{m+Q_{(0)}}\lambda_{(0,i)}^2\rt\\
    &=2m(m-1)(2m-1)
    +\frac{3m^2((1+\alpha)m+(1-\alpha))}{\alpha}\\
    &\quad+\frac{3m((1+\alpha)m+(1-\alpha))((1+\alpha)m+2(1-\alpha))}
    {\alpha(3\alpha-1)}.
\end{align*}
Finally, $\bbE_\mu[N(N-1)]=\bbE_\mu[M]$. Thus, by Corollary~\ref{cor:cc_global}(c) and \eqref{eq:cor_cc_global_tree},
\[
    \cc_{G_n}\convp\cc_{G^*}
    =\frac{3\bbE_\mu[M]p+\bbE_\mu[P]p^3}
    {\bbE_\mu[M]+2\bbE_\mu[MN]p+(\bbE_\mu[M^2]-\bbE_\mu[M])p^2}.
\]
More precisely,
\begin{equation}\label{eq:cc_PA}
    \cc_{G_n}\convp \frac{t_{\alpha, p}}{d_{\alpha, p}},
\end{equation}
where
\begin{align}
    t_{\alpha,p}
    &:=12\alpha^2(3\alpha-1)\lf(4\alpha+1)m+1-2\alpha\rf p\nonumber\\
    &\quad+8\alpha^2\lf4\alpha(3\alpha+1)m^2
    +3(1+3\alpha-6\alpha^2)m+6\alpha^2-7\alpha+3\rf p^3, \label{eq:cc_PA_t}\\
    d_{\alpha,p}
    &:=4\alpha^2(3\alpha-1)\lf(4\alpha+1)m+1-2\alpha\rf \nonumber\\
    &\quad+8\alpha(3\alpha-1)
    \lf(8\alpha^2+5\alpha+1)m^2
    +(1+2\alpha-4\alpha^2)m-\alpha\rf p\nonumber\\
    &\quad+\bigg( (3\alpha-1)(64\alpha^3+69\alpha^2+30\alpha+5)m^3\nonumber\\
    &\qquad+(-192\alpha^4+35\alpha^3+61\alpha^2+5\alpha-5)m^2\nonumber\\
    &\qquad+(1+11\alpha+7\alpha^2-51\alpha^3)m+24\alpha^4-15\alpha^3+7\alpha^2-\alpha+1\bigg) p^2 \label{eq:cc_PA_d}. 
\end{align}
In Figure~\ref{fig:pa_global_clustering_surfaces}, we plot $\cc_{G^*}$ as a function of $m$ and $p$ for $\alpha=0.4$ and $0.8$. Note that our argument does not provide convergence of $\cc_{G_n}$ for $\alpha = 0.4$.

\begin{remark}\label{rmk:alpha_condition}
    The condition $\alpha > 1/2$ is used to verify the requirement of Remark~\ref{rmk:4_moment_UI} which implies the uniform integrability assumption of Corollary~\ref{cor:cc_global}, whereas the limiting expression is finite for $\alpha > 1/3$. Thus, our result seems suboptimal in this case.
\end{remark}

\section*{Declaration of LLM assistance}
Gemini~3.1 Pro was used for the following tasks. It provided relevant references which helped the author write Section~\ref{sec:intro}. It also suggested the dataset used in Section~\ref{sec:empirical}.  The Python code required for the data analysis in Section~\ref{sec:empirical} was initially written by the LLM. Then, it was edited by the author. Gemini~3.1 Pro was further used to check the algebraic calculations in subsection~\ref{subsec:PA} which was initially done by the author. The model construction, theorem statements and proofs are solely done by the author without any help from LLM's.

\section*{Acknowledgements}
The author thanks Prof. Shirshendu Chatterjee, Prof. Remco van der Hofstad and Dr. Soumendu Sundar Mukherjee for helpful discussions about the problem. The author is also grateful to International Centre for Theoretical Sciences, Bengaluru where the problem first originated through discussions with other researchers during the ICTS-Networks workshop ``Challenges in Networks" (2026).

\bibliographystyle{alpha}
\bibliography{refs}

\appendix
\section{Proof of the supporting lemmas}\label{sec:aux}

\begin{proof}[Proof of Lemma~\ref{lem:conditional_expectation_bernoulli}]
    For $0 \le s \le m$, by independence and \eqref{eq:J_defn},
    \[
        \bbP(\Xi=s) = \sum_{\substack{A \subseteq [m]\\|A|=s}} \prod_{k \in A} r_k \prod_{k \in [m] \setminus A} (1-r_k) = J_s\lf \lf r_k \rf_{k \in [m]} \rf.
    \]
    Now, by independence,
    \begin{align*}
        \bbP(\xi_i=1, \Xi=s) &= r_i \bbP\lf \sum_{k \in [m] \setminus \{i\}} \xi_k=s-1 \rf = r_i J_{s-1}\lf \lf r_k \rf_{k \in [m] \setminus \{i\}} \rf,\\
        \bbP(\xi_i=\xi_j=1, \Xi=s) &= r_i r_j \bbP\lf \sum_{k \in [m] \setminus \{i,j\}} \xi_k=s-2 \rf = r_i r_j J_{s-2}\lf \lf r_k \rf_{k \in [m] \setminus \{i,j\}} \rf.
    \end{align*}
    Dividing by $\bbP(\Xi=s)$ gives \eqref{eq:conditional_xi_1} and \eqref{eq:conditional_xi_1_xi_2} whenever $\bbP(\Xi=s)>0$. Lastly, if $r_i$'s are all equal, \eqref{eq:conditional_xi_1_equal} and \eqref{eq:conditional_xi_1_xi_2_equal} follow by exchangeability, since $\sum_i \xi_i=\Xi$ and $\sum_{i \ne j} \xi_i \xi_j=\Xi(\Xi-1)$.
\end{proof}

We now prove Lemma~\ref{lem:expectation_computations_PA}. We shall require the so called \emph{Campbell's formula} in the proofs of \eqref{eq:estimate_6}, \eqref{eq:estimate_7}, \eqref{eq:estimate_8} and \eqref{eq:estimate_9}. We mention it here. The reader is referred to \cite{chiu2013stochastic} for details.
\begin{theorem}[Campbell's formula]\label{thm:campbell}
    Let $d,q\in\bbN$. Suppose $\Lambda$ is an inhomogeneous Poisson point process on $\bbR^d$ with locally integrable intensity function $\lambda:\bbR^d\to[0,\infty)$. Let $f:(\bbR^d)^q\to[0,\infty]$ be measurable. Then,
    \[
        \bbE\left[
          \sum_{\substack{x_1,\ldots,x_q\in\Lambda\\x_i\text{'s} \,\, \text{pairwise distinct}}}
                 f(x_1,\ldots,x_q)\right]
        =\int_{(\bbR^d)^q}f(y_1,\ldots,y_q)
                 \prod_{j=1}^q\lambda(y_j)\,dy_1\cdots dy_q.
    \]
    The sum is over ordered tuples. The identity also holds for signed $f$ when the integral of $|f|$ against the displayed product intensity is finite.
\end{theorem}
\begin{proof}[Proof of Lemma~\ref{lem:expectation_computations_PA}]
\textbf{Proof of \eqref{eq:estimate_1}:}
\begin{align*}
    \bbE_\mu[\lambda_{(0)}] &= \bbE_\mu[T_{(0)}] \left( \bbE_\mu \left[ U_{(0)}^{-\frac{1-\alpha}{2}} \right] -1 \right) \\
    &= \frac{(1 + \alpha)m}{1- \alpha} \left( \frac{2}{1 + \alpha} -1 \right)\\
    &= m.
\end{align*}
\textbf{Proof of \eqref{eq:estimate_2}:}
\begin{align*}
    \bbE_\mu[\lambda_{(0)}^2] &= \bbE_\mu \left[ T_{(0)}^2 \right] \bbE_\mu \left[ \left( x_{(0)}^{-\frac{1-\alpha}{1+ \alpha}} - 1 \right)^2 \right] \\
    &= \bbE_\mu \left[ T_{(0)}^2 \right] \bbE_\mu \left[ \left( U_{(0)}^{-\frac{1-\alpha}{2}} - 1 \right)^2 \right] \\
    &= \left( \frac{(1 + \alpha)m}{1- \alpha} + \frac{(1 + \alpha)^2 m^2}{(1- \alpha)^2} \right) \left( \frac{1}{\alpha} +1 - \frac{4}{1 + \alpha} \right) \\
    &= \frac{m ((1 + \alpha) m + (1 - \alpha))}{ \alpha}.
\end{align*}
\textbf{Proof of \eqref{eq:estimate_3}:}
\begin{align*}
    \bbE_\mu \lt \lambda_{(0,1)} \rt 
    &= \bbE_\mu \lt T_{(0,1)} \rt \lf \bbE_\mu \lt x_{(0,1)}^{-\frac{1- \alpha}{1 + \alpha}} \rt - 1 \rf \\
    &= \lf \frac{(1 + \alpha)m}{1 - \alpha} +1 \rf \lf \bbE_\mu \lt U_{(0)}^{- \frac{1 - \alpha}{2}} U_{(0,1)}^{-\frac{1- \alpha}{1 + \alpha}} \rt - 1 \rf \\
    &= \lf \frac{(1 + \alpha)m}{1 - \alpha} +1 \rf \lf \frac{1}{\alpha} - 1 \rf \\
    &= \frac{(1 + \alpha) m + (1 - \alpha)}{ \alpha}.
\end{align*}
\textbf{Proof of \eqref{eq:estimate_4}:}
\begin{align*}
    \bbE_\mu \lt \lambda_{(0)} \lambda_{(0,1)} \rt 
    &= \bbE_\mu \lt T_{(0)} \left( x_{(0)}^{-\frac{1 - \alpha}{1 + \alpha}} - 1 \right) T_{(0,1)} \left( x_{(0,1)}^{-\frac{1- \alpha}{1 + \alpha}} -1 \right) \rt \\
    &= \frac{(1 + \alpha)m}{1 - \alpha} \lf \frac{(1 + \alpha)m}{1 - \alpha} + 1 \rf \bbE_\mu \lt \lf U_{(0)}^{- \frac{1 - \alpha}{2}} - 1 \rf \lf U_{(0)}^{-\frac{1 - \alpha}{2}} U_{(0,1)}^{- \frac{1 - \alpha}{1 + \alpha}} - 1 \rf \rt \\
    &= \frac{(1 + \alpha)m}{1 - \alpha} \lf \frac{(1 + \alpha)m}{1 - \alpha} + 1 \rf \frac{(1 - \alpha)^2 (1 + 2 \alpha)}{2 \alpha^2 (1 + \alpha)} \\
    &= \frac{m (1 + 2 \alpha) ((1 + \alpha)m + (1 - \alpha))}{2 \alpha^2}.
\end{align*}
\textbf{Proof of \eqref{eq:estimate_4.5}:}
\begin{align*}
    &\bbE_\mu \lt \lambda_{(0,1)}^2 \rt \\
    &= \bbE_\mu \lt T_{(0,1)}^2 \rt \bbE_\mu \lt \lf x_{(0,1)}^{-\frac{1- \alpha}{1 + \alpha}} - 1 \rf^2 \rt \\
    &= \lf \frac{(1 + \alpha)m}{1 - \alpha} +1 \rf \lf \frac{(1 + \alpha)m}{1 - \alpha} + 2 \rf \bbE_\mu \lt \lf U_{(0)}^{- \frac{1 - \alpha}{2}} U_{(0,1)}^{-\frac{1- \alpha}{1 + \alpha}} - 1 \rf ^2 \rt \\
    &= \lf \frac{(1 + \alpha)m}{1 - \alpha} +1 \rf \lf \frac{(1 + \alpha)m}{1 - \alpha} + 2 \rf \bbE_\mu \lt U_{(0)}^{- (1 - \alpha)} U_{(0,1)}^{-\frac{2(1- \alpha)}{1 + \alpha}} - 2 U_{(0)}^{- \frac{1 - \alpha}{2}} U_{(0,1)}^{-\frac{1- \alpha}{1 + \alpha}} + 1 \rt \\
    &= \lf \frac{(1 + \alpha)m}{1 - \alpha} +1 \rf \lf \frac{(1 + \alpha)m}{1 - \alpha} + 2 \rf \lf \frac{1}{\alpha} \frac{1 + \alpha}{- 1 + 3 \alpha} - 2 \frac{2}{1 + \alpha} \frac{1 + \alpha}{2 \alpha} + 1 \rf \qquad \lt \text{since}\,\, \alpha > \frac{1}{3} \rt \\
    &= \frac{3 ((1 + \alpha)m + (1 - \alpha)) ((1 + \alpha)m + 2 (1 - \alpha))}{\alpha (-1 + 3 \alpha)}.
\end{align*}
\textbf{Proof of \eqref{eq:estimate_5}:}

\begin{align*}
    &\bbE_\mu \lt \lambda_{(0,1)} \lambda_{(0,2)} \rt \\
    &= \bbE_\mu \lt T_{(0,1)} \left( x_{(0,1)}^{-\frac{1- \alpha}{1 + \alpha}} -1 \right) T_{(0,2)} \left( x_{(0,2)}^{-\frac{1- \alpha}{1 + \alpha}} -1 \right) \rt \\
    &= \lf \frac{(1 + \alpha)m}{1 - \alpha} +1 \rf^2 \bbE_\mu \lt U_{(0)}^{- (1 - \alpha)} U_{(0,1)}^{-\frac{1- \alpha}{1 + \alpha}} U_{(0,2)}^{-\frac{1- \alpha}{1 + \alpha}} - U_{(0)}^{- \frac{1 - \alpha}{2}} U_{(0,1)}^{-\frac{1- \alpha}{1 + \alpha}} - U_{(0)}^{- \frac{1 - \alpha}{2}} U_{(0,2)}^{-\frac{1- \alpha}{1 + \alpha}} + 1 \rt \\
    &= \lf \frac{(1 + \alpha)m}{1 - \alpha} +1 \rf^2 \lf \frac{1}{\alpha} \lf \frac{1 + \alpha}{2 \alpha} \rf^2 - 2 \frac{2}{1 + \alpha} \frac{1 + \alpha}{2 \alpha} + 1 \rf \\
    &= \frac{(1 + 4 \alpha)((1 + \alpha)m + (1 - \alpha))^2}{4 \alpha^3} ,
\end{align*}
Consequently, with the mixed term omitted when $m=1$,
\begin{align*}
    &\bbE_\mu \lt \lf  \sum_{i=1}^m Q_{(0,i)}\rf^2 \rt \\
    &= \bbE_\mu \lt \lf \sum_{i=1}^m \lambda_{(0,i)} \rf^2 + \lf \sum_{i=1}^m \lambda_{(0,i)} \rf \rt \\
    &= m \bbE_\mu \lt \lambda_{(0,1)} \rt + m \bbE_\mu [\lambda_{(0,1)}^2] + m (m-1) \bbE_\mu \lt \lambda_{(0,1)} \lambda_{(0,2)} \rt \\
    &= m \frac{(1 + \alpha) m + (1 - \alpha)}{ \alpha} + m \frac{3 ((1 + \alpha)m + (1 - \alpha)) ((1 + \alpha)m + 2 (1 - \alpha))}{\alpha (-1 + 3 \alpha)} \\
    & + m (m-1) \frac{(1 + 4 \alpha)((1 + \alpha)m + (1 - \alpha))^2}{4 \alpha^3} \\
    &= \frac{m ((1 + \alpha)m + (1 - \alpha)) ((-1 - 2 \alpha + 11 \alpha^2 + 12 \alpha^3) m^2 + (2 \alpha + 14 \alpha^2 - 12 \alpha^3) m + (1 + 7 \alpha^2))}{4 \alpha^3 (-1 + 3 \alpha)}.
\end{align*}
\textbf{Proof of \eqref{eq:estimate_6}:} 
\begin{align*}
    \bbE_\mu \lt \sum_{i=m+1}^{m+ Q_{(0)}} \lambda_{(0,i)} \rt 
    &= \bbE_\mu \lt \sum_{i=m+1}^{m+ Q_{(0)}} T_{(0,i)} \left( x_{(0,i)}^{-\frac{1- \alpha}{1 + \alpha}} -1 \right) \rt \\
    &= \frac{(1 + \alpha)m}{1- \alpha} \bbE_\mu \lt \sum_{i=m+1}^{m+ Q_{(0)}} x_{(0,i)}^{-\frac{1- \alpha}{1 + \alpha}} - Q_{(0)} \rt \\
    &= \frac{(1 + \alpha)m}{1- \alpha} \bbE_\mu \lt \int_{x_{(0)}}^1 x^{- \frac{1 - \alpha}{1 + \alpha}} \rho_{(0)}(x) \, dx - \lambda_{(0)} \rt \\
    &= \frac{(1 + \alpha)m}{1- \alpha} \bbE_\mu \lt \frac{1 - \alpha}{1 + \alpha} T_{(0)} x_{(0)}^{-\frac{1 - \alpha}{1 + \alpha}}  \int_{x_{(0)}}^1 x^{- \frac{1 - \alpha}{1 + \alpha}} x^{-\frac{2 \alpha}{1 + \alpha}} \, dx - \lambda_{(0)} \rt \\
    &= \frac{(1 + \alpha)m}{1- \alpha} \bbE_\mu \lt \frac{1 - \alpha}{1 + \alpha} T_{(0)} x_{(0)}^{-\frac{1 - \alpha}{1 + \alpha}} \lf - \log x_{(0)} \rf - \lambda_{(0)} \rt \\
    &= \frac{(1 + \alpha)m}{1- \alpha} \lf - \frac{1 - \alpha}{1 + \alpha} \bbE_\mu \lt T_{(0)} \rt \frac{1 + \alpha}{2} \bbE_\mu \lt U_{(0)}^{-\frac{1 - \alpha}{2}} \log U_{(0)} \rt - \bbE_\mu \lt \lambda_{(0)} \rt \rf \\
    &= \frac{(1 + \alpha)m}{1- \alpha} \lf \frac{1 - \alpha}{2} \frac{(1 + \alpha)m}{1 - \alpha} \frac{4}{(1 + \alpha)^2} - m \rf \\
    &= m^2.
\end{align*}
\textbf{Proof of \eqref{eq:estimate_7}:}
\begin{align*}
    &\bbE_\mu \lt \lambda_{(0,1)} \sum_{i= m +1}^{m + Q_{(0)}} \lambda_{(0,i)} \rt \\
    &= \bbE_\mu \lt T_{(0,1)} \left( x_{(0,1)}^{-\frac{1- \alpha}{1 + \alpha}} -1 \right) \sum_{i=m+1}^{m+ Q_{(0)}} T_{(0,i)} \left( x_{(0,i)}^{-\frac{1- \alpha}{1 + \alpha}} -1 \right) \rt \\
    &= \bbE_\mu \lt T_{(0,1)} \rt \bbE_\mu \lt T_{(0, m+1)} \rt \bbE_\mu \lt \left( x_{(0,1)}^{-\frac{1- \alpha}{1 + \alpha}} -1 \right) \lf \sum_{i=m+1}^{m+ Q_{(0)}} x_{(0,i)}^{-\frac{1- \alpha}{1 + \alpha}} - Q_{(0)} \rf \rt \\
    &= \lf \frac{(1 + \alpha)m}{1 - \alpha} + 1 \rf \frac{(1 + \alpha)m}{1- \alpha} \\
    &\qquad\times \bbE_\mu \lt \lf x_{(0)}^{-\frac{1 - \alpha}{1 + \alpha}} \frac{1 + \alpha}{2 \alpha} - 1 \rf \lf \int_{x_{(0)}}^1 x^{- \frac{1 - \alpha}{1 + \alpha}} \rho_{(0)}(x) \, dx - \lambda_{(0)} \rf \rt  \\
    &= \frac{m(1 + \alpha)((1 + \alpha)m + (1 - \alpha))}{(1 - \alpha)^2} \\
    &\qquad\times \bbE_\mu \lt \lf x_{(0)}^{-\frac{1 - \alpha}{1 + \alpha}} \frac{1 + \alpha}{2 \alpha} - 1 \rf \lf \frac{1 - \alpha}{1 + \alpha} T_{(0)} x_{(0)}^{-\frac{1 - \alpha}{1 + \alpha}} \lf - \log x_{(0)} \rf - \lambda_{(0)} \rf \rt \\
    &= \frac{m(1 + \alpha)((1 + \alpha)m + (1 - \alpha))}{(1 - \alpha)^2} \bbE_\mu \bigg[ \frac{1 - \alpha}{2 \alpha} T_{(0)} x_{(0)}^{-\frac{2(1 - \alpha)}{1 + \alpha}} (- \log x_{(0)}) \\
    &- \frac{1 + \alpha}{2 \alpha} \lambda_{(0)} x_{(0)}^{-\frac{1 - \alpha}{1 + \alpha}} - \frac{1 - \alpha}{1 + \alpha} T_{(0)} x_{(0)}^{-\frac{1 - \alpha}{1 + \alpha}} \lf - \log x_{(0)} \rf + \lambda_{(0)} \bigg] \\
    &=\frac{m(1+\alpha)((1+\alpha)m+(1-\alpha))}{(1-\alpha)^2}
    \left(\frac{m(1+\alpha)^2}{4\alpha^3}
    -\frac{m(1+\alpha)}{2\alpha^2}-\frac{2m}{1+\alpha}+m\right) \\
    &= \frac{m(1 + \alpha)((1 + \alpha)m + (1 - \alpha))}{(1 - \alpha)^2} \frac{m (1 - \alpha)^2 (4 \alpha^2 + 3 \alpha + 1)}{4 \alpha^3 (1 + \alpha)} \\
    &= \frac{m^2 ((1 + \alpha) m + (1 - \alpha)) (1 + 3 \alpha + 4 \alpha^2)}{4 \alpha^3}.
\end{align*}
\textbf{Proof of \eqref{eq:estimate_8}:}
\begin{align*}
    &\bbE_\mu \lt Q_{(0)} \sum_{i= m +1}^{m + Q_{(0)}} \lambda_{(0,i)} \rt \\
    &= \bbE_\mu \lt \sum_{i, j = m +1}^{m + Q_{(0)}} \lambda_{(0,i)} \rt \\
    &= \bbE_\mu \lt \sum_{i = m +1}^{m + Q_{(0)}} \lambda_{(0,i)} \rt + \bbE_\mu \lt \sum_{\substack{ i, j = m +1\\ i \ne j }}^{m + Q_{(0)}} \lambda_{(0,i)} \rt \\
    &= m^2 + \bbE_\mu \lt \sum_{\substack{ i, j =m+1 \\ i \ne j }}^{m+ Q_{(0)}} T_{(0,i)} \left( x_{(0,i)}^{-\frac{1- \alpha}{1 + \alpha}} -1 \right) \rt \\
    &= m^2 + \bbE_\mu \lt T_{(0, m+1)} \rt \bbE_\mu \lt \sum_{\substack{ i, j =m+1 \\ i \ne j }}^{m+ Q_{(0)}} x_{(0,i)}^{-\frac{1- \alpha}{1 + \alpha}} - Q_{(0)} \lf Q_{(0)} - 1 \rf \rt \\
    &= m^2 + \frac{(1 + \alpha) m}{1 - \alpha} \bbE_\mu \lt \int_{x_{(0)}}^1 \int_{x_{(0)}}^1 x^{-\frac{1 - \alpha}{1 + \alpha}} \rho_{(0)}(x) \rho_{(0)}(y) \, dx \, dy - \lambda_{(0)}^2 \rt  \\
    &= m^2 + \frac{(1 + \alpha) m}{1 - \alpha} \bbE_\mu \lt \frac{1 - \alpha}{1 + \alpha} T_{(0)} x_{(0)}^{-\frac{1 - \alpha}{1 + \alpha}} \lf - \log x_{(0)} \rf \lambda_{(0)} - \lambda_{(0)}^2 \rt \\
    &= m^2 + \frac{(1 + \alpha) m}{1 - \alpha} \lf \frac{1 - \alpha}{1 + \alpha} \bbE_\mu \lt T_{(0)}^2 \rt \bbE_\mu \lt x_{(0)}^{-\frac{1 - \alpha}{1 + \alpha}} \lf x_{(0)}^{-\frac{1 - \alpha}{1 + \alpha}} - 1 \rf \lf - \log x_{(0)} \rf \rt - \bbE_\mu \lt \lambda_{(0)}^2 \rt \rf \\
    &= m^2 + \frac{(1 + \alpha) m}{1 - \alpha} \lf \frac{1 - \alpha}{1 + \alpha} \frac{(1 + \alpha) m}{1 - \alpha} \lf \frac{(1 + \alpha) m}{1 - \alpha} + 1 \rf \frac{(1 + 3 \alpha) (1 - \alpha)}{2 \alpha^2 (1 + \alpha)} - \frac{m ((1 + \alpha) m + (1 - \alpha))}{ \alpha} \rf \\
    &= \frac{m^2 (1 + \alpha) ( (1 + 2 \alpha) m +1)}{2 \alpha^2}.
\end{align*}
\textbf{Proof of \eqref{eq:estimate_8.5}:}
\begin{align*}
    &\bbE_\mu \lt \sum_{i=m+1}^{m+ Q_{(0)}} \lambda_{(0,i)}^2 \rt \\
    &= \bbE_\mu \lt \sum_{i=m+1}^{m+ Q_{(0)}} T_{(0,i)}^2 \left( x_{(0,i)}^{-\frac{1- \alpha}{1 + \alpha}} -1 \right)^2 \rt \\
    &= \frac{(1 + \alpha)m}{1- \alpha} \lf \frac{(1 + \alpha)m}{1- \alpha} + 1 \rf \bbE_\mu \lt \sum_{i=m+1}^{m+ Q_{(0)}} x_{(0,i)}^{-\frac{2(1- \alpha)}{1 + \alpha}} - 2 \sum_{i=m+1}^{m+ Q_{(0)}} x_{(0,i)}^{-\frac{1- \alpha}{1 + \alpha}} + Q_{(0)} \rt \\
    &= \frac{(1 + \alpha)m}{1- \alpha} \lf \frac{(1 + \alpha)m}{1- \alpha} + 1 \rf \bbE_\mu \lt \int_{x_{(0)}}^1 x^{- \frac{2 (1 - \alpha)}{1 + \alpha}} \rho_{(0)}(x) \, dx -2 \int_{x_{(0)}}^1 x^{- \frac{1 - \alpha}{1 + \alpha}} \rho_{(0)}(x) \, dx + \lambda_{(0)} \rt \\
    &= \frac{(1 + \alpha)m}{1- \alpha} \lf \frac{(1 + \alpha)m}{1- \alpha} + 1 \rf \bbE_\mu \bigg[ \frac{1 - \alpha}{1 + \alpha} T_{(0)} x_{(0)}^{-\frac{1 - \alpha}{1 + \alpha}}  \int_{x_{(0)}}^1 x^{- \frac{2(1 - \alpha)}{1 + \alpha}} x^{-\frac{2 \alpha}{1 + \alpha}} \, dx \\
    & - \frac{2(1 - \alpha)}{1 + \alpha} T_{(0)} x_{(0)}^{-\frac{1 - \alpha}{1 + \alpha}}  \int_{x_{(0)}}^1 x^{- \frac{1 - \alpha}{1 + \alpha}} x^{-\frac{2 \alpha}{1 + \alpha}} \, dx + \lambda_{(0)} \bigg] \\
    &= \frac{(1 + \alpha)m}{1- \alpha} \lf \frac{(1 + \alpha)m}{1- \alpha} + 1 \rf \bbE_\mu \bigg[ \frac{1 - \alpha}{1 + \alpha} T_{(0)} x_{(0)}^{-\frac{1 - \alpha}{1 + \alpha}} \frac{1 + \alpha}{1 - \alpha} \lf x_{(0)}^{ - \frac{1 - \alpha}{1 + \alpha}} - 1 \rf \\
    & - \frac{2(1 - \alpha)}{1 + \alpha} T_{(0)} x_{(0)}^{-\frac{1 - \alpha}{1 + \alpha}} \lf - \log x_{(0)} \rf + \lambda_{(0)} \bigg] \\
    &= \frac{(1 + \alpha)m}{1- \alpha} \lf \frac{(1 + \alpha)m}{1- \alpha} + 1 \rf \bigg( \bbE_\mu \lt T_{(0)} \rt \lf \bbE_\mu \lt U_{(0)}^{-(1 - \alpha)} \rt - \bbE_\mu \lt U_{(0)}^{-\frac{1-\alpha}{2}} \rt \rf \\
    & - \frac{2 (1 - \alpha)}{1 + \alpha} \bbE_\mu \lt T_{(0)} \rt \frac{1 + \alpha}{2} \bbE_\mu \lt U_{(0)}^{-\frac{1 - \alpha}{2}} ( - \log U_{(0)}) \rt + \bbE_\mu \lt \lambda_{(0)} \rt \bigg) \\
    &= \frac{(1 + \alpha)m}{1- \alpha} \lf \frac{(1 + \alpha)m}{1- \alpha} + 1 \rf \bigg( \frac{(1 + \alpha) m}{1 - \alpha} \lf \frac{1}{\alpha} - \frac{2}{1 + \alpha} \rf \\
    & - (1 - \alpha) \frac{(1 + \alpha)m}{1 - \alpha} \frac{4}{(1 + \alpha)^2} + m \bigg) \\ 
    &= \frac{(1 + \alpha)m}{1- \alpha} \lf \frac{(1 + \alpha)m}{1- \alpha} + 1 \rf \frac{(1 - \alpha)^2 m}{\alpha (1 + \alpha)} \\
    &= \frac{m^2 ((1 + \alpha) m + (1 - \alpha))}{\alpha}.
\end{align*}
\textbf{Proof of \eqref{eq:estimate_9}:}
\begin{align*}
    & \bbE_\mu \lt \sum_{\substack{i, j= m +1 \\ i \ne j}}^{m + Q_{(0)}} \lambda_{(0,i)} \lambda_{(0,j)} \rt \\
    &= \bbE_\mu \lt \sum_{\substack {i, j=m+1 \\ i \ne j}}^{m+ Q_{(0)}} T_{(0,i)} T_{(0,j)} \left( x_{(0,i)}^{-\frac{1- \alpha}{1 + \alpha}} -1 \right) \left( x_{(0,j)}^{-\frac{1- \alpha}{1 + \alpha}} -1 \right) \rt \\
    &= \lf \frac{(1 + \alpha) m}{1 - \alpha} \rf^2 \bbE_\mu \lt \int_{x_{(0)}}^1 \int_{x_{(0)}}^1 \lf x^{-\frac{1 - \alpha}{1 + \alpha}} - 1 \rf \lf y^{-\frac{1 - \alpha}{1 + \alpha}} - 1 \rf \rho_{(0)}(x) \rho_{(0)} (y) \, dx \, dy \rt \\
    &= \lf \frac{(1 + \alpha) m}{1 - \alpha} \rf^2 \bbE_\mu \lt \lf \int_{x_{(0)}}^1 \lf x^{-\frac{1 - \alpha}{1 + \alpha}} - 1 \rf \rho_{(0)}(x) \, dx \rf^2 \rt \\
    &= \lf \frac{(1 + \alpha) m}{1 - \alpha} \rf^2 \bbE_\mu \lt \lf \frac{1 - \alpha}{1 + \alpha} T_{(0)} x_{(0)}^{-\frac{1 - \alpha}{1 + \alpha}} \lf - \log x_{(0)} \rf - \lambda_{(0)}\rf^2 \rt \\
    &= \lf \frac{(1 + \alpha) m}{1 - \alpha} \rf^2 \bbE_\mu \bigg[ \lf \frac{1 - \alpha}{1 + \alpha} \rf^2 T_{(0)}^2  x_{(0)}^{-\frac{2 (1 - \alpha)}{1 + \alpha}} \lf \log x_{(0)} \rf^2 \\
    &\qquad - \frac{2(1 - \alpha)}{1 + \alpha} T_{(0)} x_{(0)}^{-\frac{1 - \alpha}{1 + \alpha}} \lf - \log x_{(0)} \rf \lambda_{(0)} + \lambda_{(0)}^2 \bigg] \\
    &= \lf \frac{(1 + \alpha) m}{1 - \alpha} \rf^2 \bigg( \lf \frac{1 - \alpha}{1 + \alpha} \rf^2 \bbE_\mu \lt T_{(0)}^2 \rt \bbE_\mu \lt x_{(0)}^{-\frac{2 (1 - \alpha)}{1 + \alpha}} \lf \log x_{(0)} \rf^2 \rt \\
    & - \frac{2(1 - \alpha)}{1 + \alpha} \bbE_\mu \lt T_{(0)} x_{(0)}^{-\frac{1 - \alpha}{1 + \alpha}} \lf - \log x_{(0)} \rf T_{(0)} \lf U_{(0)}^{-\frac{1 - \alpha}{2}} - 1 \rf \rt + \bbE_\mu \lt \lambda_{(0)}^2 \rt \bigg) \\
    &= \lf \frac{(1 + \alpha) m}{1 - \alpha} \rf^2 \bigg( \lf \frac{1 - \alpha}{1 + \alpha} \rf^2 \frac{(1 + \alpha) m}{1 - \alpha} \lf \frac{(1 + \alpha) m}{1 - \alpha} + 1 \rf \bbE_\mu \lt x_{(0)}^{-\frac{2 (1 - \alpha)}{1 + \alpha}} \lf \log x_{(0)} \rf^2 \rt \\
    & - \frac{2(1 - \alpha)}{1 + \alpha} \frac{(1 + \alpha) m}{1 - \alpha} \lf \frac{(1 + \alpha) m}{1 - \alpha} + 1 \rf\bbE_\mu \lt x_{(0)}^{-\frac{1 - \alpha}{1 + \alpha}} \lf - \log x_{(0)} \rf \left( U_{(0)}^{-\frac{1 - \alpha}{2}} - 1 \right) \rt \\
    &\qquad + \frac{m ((1 + \alpha) m + (1 - \alpha))}{ \alpha} \bigg) \\
    &=\left(\frac{(1+\alpha)m}{1-\alpha}\right)^2
       \bigg(\frac{m(1+\alpha)((1+\alpha)m+(1-\alpha))}{2\alpha^3}\\
    &\qquad-\frac{m(1+3\alpha)((1+\alpha)m+(1-\alpha))}{\alpha^2(1+\alpha)}
          +\frac{m((1+\alpha)m+(1-\alpha))}{\alpha}\bigg)\\
    &=\frac{m^3(1+\alpha)(1+2\alpha)((1+\alpha)m+(1-\alpha))}{2\alpha^3},
\end{align*}
where we used
\[
    \bbE_\mu\lt x_{(0)}^{-\frac{2(1-\alpha)}{1+\alpha}} (\log x_{(0)})^2\rt
    =\frac{(1+\alpha)^2}{4}\int_0^1u^{\alpha-1}(\log u)^2\,du
    =\frac{(1+\alpha)^2}{2\alpha^3}.
\]
Finally,
\begin{align*}
    &\bbE_\mu\lt\left(\sum_{i=m+1}^{m+Q_{(0)}}Q_{(0,i)}\right)^2\rt\\
    &=\bbE_\mu\lt\sum_{i=m+1}^{m+Q_{(0)}}\lambda_{(0,i)}\rt +\bbE_\mu\lt\sum_{i=m+1}^{m+Q_{(0)}}\lambda_{(0,i)}^2\rt
    +\bbE_\mu\lt\sum_{\substack{i,j=m+1\\i\ne j}}^{m+Q_{(0)}}
    \lambda_{(0,i)}\lambda_{(0,j)}\rt\\
    &=m^2+\frac{m^2((1+\alpha)m+(1-\alpha))}{\alpha}
    +\frac{m^3(1+\alpha)(1+2\alpha)((1+\alpha)m+(1-\alpha))}{2\alpha^3}.
\end{align*}
This proves \eqref{eq:estimate_9} and completes the proof.
\end{proof}

\end{document}